\documentclass[leqno,12pt]{article} 
\usepackage{amsmath, amssymb}
\usepackage{amsthm} 
\usepackage{epsfig}
\usepackage{lipsum}
\usepackage{graphicx}
\usepackage{calc}
\usepackage{gensymb}
\usepackage{comment}
\usepackage{wrapfig}
\usepackage[sort, numbers]{natbib}
\usepackage{placeins}
\usepackage{hyperref}
\usepackage[dvipsnames]{xcolor}
\usepackage{ulem}
\usepackage{bbm}
\usepackage{authblk}
\usepackage{etex}
\usepackage[etex=true,export]{adjustbox}
\usepackage{float}
\usepackage{xkeyval}
\usepackage{mathtools}
\usepackage{amsmath}
\usepackage{amsthm}
\usepackage{amsfonts}
\usepackage{etoolbox}
\usepackage{hyperref}
\usepackage{amssymb}
\usepackage{array}
\usepackage{multirow}
\usepackage{lipsum}
\usepackage{hyperref}
\usepackage{authblk}
\usepackage{url}
\usepackage{epstopdf,epsfig}

\theoremstyle{plain} 
\newtheorem{theorem}{\indent\sc Theorem}[section]
\newtheorem{lemma}[theorem]{\indent\sc Lemma}

\newtheorem{proposition}[theorem]{\indent\sc Proposition}

\theoremstyle{definition} 
\newtheorem{definition}[theorem]{\indent\sc Definition}

\title{\uppercase{\normalsize Hamiltonian Properties of Hybrid-Faulty Burnt Pancake Graphs}} 
\author{
    Hongyi Zhu, Qingying Deng\footnote{Corresponding author. \newline {\em E-mail address:} qingying@xtu.edu.cn (Q. Deng).} \\
    {\footnotesize Key Laboratory of Intelligent Computing and Information Processing of Education,}\\
    {\footnotesize School of Mathematics and Computational Science, Xiangtan University, Xiangtan, Hunan 411105, PR China}
}

\date{}

\begin{document}
\maketitle


\begin{abstract} 
We investigate the combined occurrence of edge faults and vertex faults in the burnt pancake graph (\( BP_n \)). In this paper, we prove that \( BP_n - F \), where \( F \) includes pairs of end-vertices of matching edges and fault-tolerant edges, contains a Hamiltonian cycle when \( |F| \leq n-2 \) and a Hamiltonian path when \( |F| \leq n-3 \). This establishes that \( BP_n \) is \((n-2)\)-hybrid fault Hamiltonian and \((n-3)\)-hybrid fault Hamiltonian connected for \( n \geq 3 \). These results are demonstrated to be optimal under the given conditions, with all bounds shown to be tight.
\end{abstract}

\medskip

\section{Introduction}\label{section1}

In multi-processor systems, interconnection networks exhibit various specific characteristics, and the study of their topological structures is crucial.  
An interconnection network (or simply, network) is represented by a connected graph \( G = (V, E) \), where \( V \) represents processors and \( E \) represents communication links between processors.  
Selecting an appropriate network is a critical task in multi-processor systems, requiring the analysis of key characteristics such as symmetry, transitivity, connectivity, distance, recursive structure, and maximal fault tolerance.  
In practical applications, failures of processors or communication links are inevitable.  
A processor (or vertex) is considered fault-free if it is operational, and a communication link (or edge) is considered fault-free if both its end-vertices and the edge itself are operational.  
The \( n \)-dimensional burnt pancake graph \( BP_n \) is a popular topology for designing networks in distributed systems and large-scale parallel computing.  
It exhibits desirable topological properties, including high symmetry, vertex transitivity, regularity, recursiveness, and maximal fault tolerance.  
The pancake problem was first introduced in \cite{1}.  
In \cite{14}, pancake graphs were formally defined as Cayley graphs on the symmetric group with generating sets consisting of all prefix-reversals.  
In this paper, we focus on Hamiltonian cycles and paths in faulty burnt pancake graphs.

The cycle network and the path network, with a minimal number of links and no redundant branches, are two popular interconnection network topologies.  
These networks are not only suitable for local area networks but also enable the development of simple and efficient algorithms with low communication costs.  
If a network is embedded in a cycle, it can execute all the algorithms applicable to cycles as well. 
As a result, embedding cycles and paths has been a subject of significant interest among researchers.  
If a path \( H \) passes through every vertex of \( G \) exactly once, \( H \) is called a \emph{Hamiltonian path} \cite{24}.  
Similarly, if a cycle \( C \) passes through every vertex of \( G \) exactly once, \( C \) is called a \emph{Hamiltonian cycle} \cite{24}.  
A graph that contains a Hamiltonian cycle is called a  \emph{Hamiltonian graph}.
A graph is \emph{Hamiltonian-connected} if for every pair of vertices there is a Hamiltonian path between the two vertices.
In \cite{15}, Blanco et al. characterized all the 8-cycles in \( BP_n \) for \( n \geq 2 \), identifying them as the smallest cycles that can be embedded in \( BP_n \).
In \cite{17}, Kanevsky et al. demonstrated that all cycles of length \( l \), where \( 6 \leq l \leq n! - 2 \) or \( l = n! \), can be embedded in pancake graphs.  
They also proposed new algorithms for embedding Hamiltonian cycles and other cycle sets.

In multi-processor systems, there exist stable processors and/or links alongside faulty processors and/or links.  
When embedding cycles and paths into a fault-prone network, it is essential for the network to avoid these faulty processors and/or links.  
Let \( F^{mv} \) be a set of pairs of end-vertices of matching edges with size at most \( m \), and let \( F^e \) be a set of faulty edges of \(G-F^{mv}\) with size at most \( k - m \), where \( m \) and \( k \) are positive integers.  
A graph \( G = (V, E) \) is called \emph{\( k \)-hybrid fault Hamiltonian} if \( G - F^{mv} - F^e \) is Hamiltonian for any \( F^{mv} \) consisting of \( m \) pairs of end-vertices of matching edges and any \( F^e \) consisting of \( k - m \) faulty edges.  
Similarly, \( G \) is called \emph{\( k \)-hybrid fault Hamiltonian connected} if \( G - F^{mv} - F^e \) is Hamiltonian connected for any \( F^{mv} \) and \( F^e \) defined as above.  
An edge or vertex is considered \emph{fault-free} if it is not broken.  
In recent years, Hung et al. \cite{2} showed that the hypercube \( Q_n \) contains a Hamiltonian path between opposite partite sets with faulty matching.  
In \cite{4}, Sun et al. extended this result by proving that \( Q_n \) contains a Hamiltonian path between opposite partite sets with faulty matching and edges.  
In \cite{3}, Yang et al. demonstrated that star graph \( S_n \) contains both a Hamiltonian cycle and a Hamiltonian path between opposite partite sets with faulty matching.  
Further, in \cite{5}, Lu and Xue generalized these results to show that \( S_n \) contains a Hamiltonian cycle and a Hamiltonian path between opposite partite sets that pass through prescribed edges with faulty matching and edges.  
In \cite{16}, Hung et al. established that pancake graphs \( P_n \) are Hamiltonian and Hamiltonian connected even in the presence of faulty vertices and/or edges.  
In this paper, we focus on the problem of embedding Hamiltonian cycles and Hamiltonian paths in burnt pancake graphs with faulty pairs of end-vertices of matching and faulty edges.  
The main result is as follows:  
\( BP_n \) is \((n-2)\)-hybrid fault Hamiltonian and \((n-3)\)-hybrid fault Hamiltonian connected.

The rest of this paper is organized as follows:
In Section \ref{section2}, we introduce notations and discuss the basic properties of the burnt pancake graph $BP_{n}$.
In Section \ref{section3}, the main results on the hybrid fault tolerance of 
$BP_{n}$ are presented and proven.
In Section \ref{section4}, we conclude the paper. 

\section{Preliminary}\label{section2}

\subsection{Terminologies and notation}\label{section2.1}

Let \( G = (V, E) \) be a given graph.  
For an edge \( e = (u, v) \), the vertices \( u \) and \( v \) are called the \emph{end-vertices} of \( e \) and \(\{u,v\}\) is called \emph{a pair of end-vertices} of \( e \).  
We say \( u \) is \emph{adjacent} to \( v \), and vice versa.  
For any two vertices \( u,v \in V(G) \), we use \( d_G(u,v) \) to denote the \emph{distance} of \( u \) and \( v \) in \( G \). 
Assume that \( |*| \) represents the cardinality of the set \( * \).  
We also use \( P[u, v] \) to represent a path between any two distinct vertices \( u \) and \( v \) in \( G \), denoted by \( \langle u, \cdots, v \rangle \).  
The notations used in this paper are summarized in Table \(\ref{B1}\).

\begin{table}[h]
  \setlength{\abovecaptionskip}{0cm}
  \setlength{\belowcaptionskip}{0.2cm}
  \caption{Some notations used in this paper}\label{B1}
  \centering
  \begin{tabular}{p{7cm}<{\centering}p{9cm}<{\raggedright}}
  \hline
  Notations & Significations \\
  \hline
  $\langle n\rangle$ & $\{1,2,\cdots,n\}$ \\
  $[n]$ & $\{1,2,\cdots,n\}\cup \{\bar{1},\bar{2},\cdots,\bar{n}\}$ \\
  $(*)_{n}$ &  the $n$-th bit of a vertex $*$ \\
  $M$ &  a matching edge set of $G$ with size $m$ \\
  $F^{mv}$ & a set contains pairs of end-vertices of matching edges, i.e., $\{\{a_{i},b_{i}\}|(a_{i},b_{i})\in M$ for $1\leq i\leq m\}$ \\
  $F^{e}$ & a set contains edges of $G-F^{mv}$ with size $k-m$\\
  $F$ & a set contains $k$ hybrid fault elements, i.e., $F=F^{mv}\cup F^{e}$ \\
  $V(F^{mv})$ & $\{a_{1},b_{1},\cdots,a_{m},b_{m}\}$ \\
  $V(F)$  &  $V(F)=V(F^{mv})\cup V(F^{e})$ \\
  $|F|$  &  $|F|=|F^{mv}|+|F^{e}|$ \\
  $BP_{n}$ & The $n$-dimensional burnt pancake graph\\
  $BP_{n}^{i}$ & the subgraph induced by $\{u_{1}u_{2}\cdots u_{n}\in V(BP_{n})|u_{n}=i\}$ for $i\in [n]$\\
  \hline
  \end{tabular}
\end{table}

In this paper, we primarily investigate the Hamiltonian properties of graphs when subjected to hybrid faulty elements.  

We use \( F^{mv} \) and \( F^e \) to represent pairs of end-vertices of matching edges and faulty edges, respectively, such that \( F = F^{mv} \cup F^e \).  
Symbolically, let \( F^{mv} = \{\{a_i, b_i\} \mid (a_i, b_i) \in M, 1 \leq i \leq k\} \), where \( 1 \leq k \leq m \), and \( F^e \subset E(G - F^{mv}) \), where \( M \) is a matching edge set in $G$.  

The vertex set \( V(F^{mv}) \) is defined as \( \{a_1, b_1, \dots, a_m, b_m\} \), and the total vertex set \( V(F) \) is \( V(F^{mv}) \cup V(F^e) \).  
The cardinality of \( F \), denoted \( |F| \), is the total number of elements in \( F \), i.e., \( |F| = |F^{mv}| + |F^e| \).  

{\bf Example:}
Consider \( BP_3 \) in Figure \(\ref{Fig8}\). Assume \( F = \{\{123, \bar{1}23\}, \{1\bar{2}3, \bar{3}2\bar{1}\}, (\bar{2}13, 213), \\(\bar{1}\bar{2}3, \bar{3}21), (\bar{1}2\bar{3}, \bar{2}1\bar{3})\} \).  
Then, we have:  
\[
\begin{aligned}
F^{mv} &= \{\{123, \bar{1}23\}, \{1\bar{2}3, \bar{3}2\bar{1}\}\}, \\
F^e &= \{(\bar{2}13, 213), (\bar{1}\bar{2}3, \bar{3}21), (\bar{1}2\bar{3}, \bar{2}1\bar{3})\}, \\
V(F^{mv}) &= \{123, \bar{1}23, 1\bar{2}3, \bar{3}2\bar{1}\}, \\
|F| &= |F^{mv}| + |F^e| = 2 + 3 = 5, \\
V(F) &= \{123, \bar{1}23, 1\bar{2}3, \bar{3}2\bar{1}, \bar{2}13, 213, \bar{1}\bar{2}3, \bar{3}21, \bar{1}2\bar{3}, \bar{2}1\bar{3}\}.
\end{aligned}
\]  

The red markings in Figure \(\ref{Fig8}\) indicate the elements of \( F \), and the blue vertices represent \( V(F) \).

\subsection{The $n$-dimensional burnt pancake graph}\label{section2.2}

For any positive integer \( n \), let \( \langle n \rangle = \{1, 2, \cdots, n\} \), and let \( u_1 u_2 \cdots u_{i-1} u_i u_{i+1} \cdots u_n \) be a permutation of \( \langle n \rangle \).  
We use \( [n] \) to denote \( \langle n \rangle \cup \{\bar{1}, \bar{2}, \dots, \bar{n}\} \), where \( \bar{i} = -i \) for \( i \in \langle n \rangle \) \cite{19}.  

A \emph{signed permutation} \cite{18} of \( \langle n \rangle \) is an \( n \)-permutation \( u_1 u_2 \cdots u_n \) that selects \( n \) elements from \( [n] \) such that \( |u_1| |u_2| \cdots |u_n| \) is a permutation of \( \langle n \rangle \).  
Then \( u = u_1 u_2 \cdots u_{i-1} u_i u_{i+1} \cdots u_n \) can be a vertex of the burnt pancake graph \( BP_n \), where \( u_i \in [n] \) and \( i \in \{1, 2, \dots, n\} \), and \( (u)_n \) denotes the \( n \)-th bit of the vertex \( u \).  

The \emph{\( k \)-th prefix-reversal} \cite{18} of a signed permutation \( u = u_1 u_2 \cdots u_{k-1} u_k u_{k+1} \cdots u_n \) of \( \langle n \rangle \) is given by  
\[
\bar{u}_k \bar{u}_{k-1} \cdots \bar{u}_1 u_{k+1} \cdots u_n.
\]  
The vertex \( \bar{u}_k \bar{u}_{k-1} \cdots \bar{u}_1 u_{k+1} \cdots u_n \) is the unique \emph{\( k \)-neighbor} of \( u \), denoted by \( k(u) \).  
Suppose \( u = \bar{2} 1 \bar{6} 4 \bar{5} 3 \). Then,  
\[
3(u) = 6 \bar{1} 2 4 \bar{5} 3 \quad \text{and} \quad 6(u) = \bar{3} 5 \bar{4} 6 \bar{1} 2.
\]

\begin{definition}\cite{19}\label{D2}
The \( n \)-dimensional burnt pancake graph is denoted by \( BP_n \).  
The vertex set of \( BP_n \) is defined as  
\[
V(BP_n) = \{ u = u_1 u_2 \cdots u_n \mid u_i \in [n], u_i \neq u_j \text{ for } i \neq j \},
\]  
where each vertex is a signed permutation of \( \langle n \rangle \).  

The edge set of \( BP_n \) is defined as  
\[
E(BP_n) = \{ (u, v) \mid u = u_1 u_2 \cdots u_k \cdots u_n, \, v = k(u) = \bar{u}_k \bar{u}_{k-1} \cdots \bar{u}_1 u_{k+1} \cdots u_n \},
\]  
where the edge \( (u, v) \) is called a \emph{\( k \)-dimensional edge}, and the vertex \( v \) is referred to as the \emph{\( k \)-neighbor} of \( u \).  

\end{definition}

\begin{figure}[h]
  \centering
  \includegraphics[width=1\textwidth]{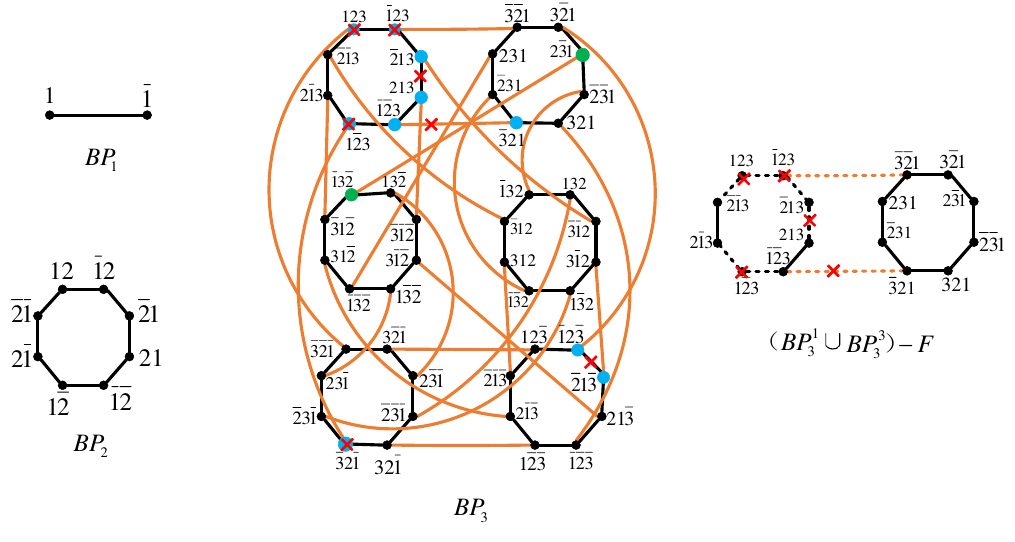}\\
  \caption{The illustrations of $BP_{1},BP_{2}$, $BP_{3}$ and $(BP_{3}^{1}\cup BP_{3}^{3})-F$.}
  \label{Fig8}
\end{figure}

First, note that \( BP_n \) is not a bipartite graph for $n\geq 3$ \cite{dvovrak2024paired}.

From the above definition, \( BP_n \) is a graph consisting of \( 2^n \cdot n! \) vertices, where each vertex represents a distinct signed permutation of \( \langle n \rangle \).  
The subgraph \( BP_n^i \), induced by the vertices \( \{ u_1 u_2 \cdots u_n \in V(BP_n) \mid u_n = i \} \), is isomorphic to \( BP_{n-1} \) for all \( i \in [n] \).  
Figure \(\ref{Fig8}\) illustrates this property by showing the subgraphs \( BP_3^1 \) and \( BP_3^3 \) of \( BP_3 \), which include fault-tolerant elements.  
The graph \( BP_n \) can be partitioned into \( 2n \) vertex-disjoint subgraphs, each isomorphic to \( BP_{n-1} \), and it is an \( n \)-regular graph with \( n \cdot n! \cdot 2^{n-1} \) edges.

\begin{lemma}\cite{20}\label{L14}
Let $n\geq 2$ and $i,j\in [n]$ such that $i\neq j$.
If $i\neq \overline{j}$, there are $(n-2)!\cdot 2^{n-2}$ edges, denoted by $E_{i,j}(BP_{n})$, between $BP_{n}^{i}$ and $BP_{n}^{j}$; Otherwise, $|E_{i,j}|=0$.
\end{lemma}

An edge of $BP_{n}$ is called an \emph{out-edge} if it is in $E_{i,j}(BP_{n})$ for $i,j\in [n]$ such that $i\neq \{j,\bar{j}\}$.
In \( BP_3 \) (as shown in Figure \(\ref{Fig8}\)), orange edges are used to denote the out-edges.

For any vertex \( u \in BP_n^i \), a neighbor \( v \) of \( u \) is called an \emph{out-neighbor} of \( u \) if \( (u, v) \) is an out-edge. 
In this paper, we classify the $n$-th bit of each vertex.
Hence, the out-neighbor of $u$ is the $n$-neighbor of $u$, and the out-neighbor \( v \) of \( u \) is unique.  
In \( BP_3 \) (as shown in Figure \(\ref{Fig8}\)), let \( s = \bar{1}3\bar{2} \in V(BP_3) \).  
Then \( 3(s) = 2\bar{3}1 \) is the out-neighbor of \( s \).  
This is illustrated by the green vertices of \( BP_3 \) in Figure \(\ref{Fig8}\).

\begin{proposition}\cite{21}\label{P2}
Let $\{k_{1}',k_{2}',\cdots,k_{m}'\}\subseteq [n]$ with $m\geq5$.
Then there is a permutation $k_{1},k_{2},\cdots,k_{m}$ of $k_{1}',k_{2}',\cdots,k_{m}'$ such that $k_{1}=k_{1}'$, $k_{m}=k_{m}'$ and $k_{i}\neq \overline{k_{i+1}} $ for $i\in \{1,2,\cdots,m-1\}$.
\end{proposition}

For example, let \( k_1' = 1 \), \( k_2' = \bar{1} \), \( k_3' = 2 \), \( k_4' = \bar{2} \), and \( k_5' = 3 \).  
Then, there exists a permutation \( k_1, k_2, \dots, k_5 \) of \( k_1', k_2', \dots, k_5' \) such that \( k_1 = 1 \), \( k_2 = 2 \), \( k_3 = \bar{1} \), \( k_4 = \bar{2} \), and \( k_5 = 3 \).  
Consequently, there are \( (n-2)! \cdot 2^{n-2} \) edges between $BP_n^{k_i}$ and $BP_n^{k_{i+1}}$ for $i=1,2,3,4$.

\begin{lemma}\cite{22}\label{L15}
Let $n\geq 3$, $u\in V(BP_{n}^{i})$ and $v\in V(BP_{n}^{j})$ for some $i,j\in [n]$.
If $i=j$ and $1\leq d_{BP_{n}}(u,v)\leq 2$, then $(n(u))_{n}\neq (n(v))_{n}$;
If $i\neq j$ and $d_{BP_{n}}(u,v)\leq 3$, then $(n(u))_{n}\neq (n(v))_{n}$.
\end{lemma}

When \( F^{mv} = \emptyset \), the following lemma is well-established.
\begin{lemma}\cite{18}\label{L16}
Let $F^{e}$ be a subset of $E(BP_{n})$ and $n$ be a positive integer with $n\geq3$.
Then $BP_{n}$ is $(n-2)$-edge fault Hamiltonian and $(n-3)$-edge fault Hamiltonian connected.
\end{lemma}
By the above lemma, we conclude that \( BP_n \) is Hamiltonian, as proven in \cite{15}.

\section{Hybrid Fault Tolerance Hamiltonian of Burnt Pancake Graphs}\label{section3}

In this section, we study the Hamiltonian properties of burnt pancake graphs with the hybrid faulty set \( F \).  
Let \( M \) be a matching edge set of \( BP_n \) with size \( m \), and let \( F^{mv} = \{\{a_i, b_i\} \mid (a_i, b_i) \in M, \ 1 \leq i \leq m\} \).  
Without loss of generality (w.l.o.g.), assume that \( F = F^{mv} \cup F^e \), where \( F^{mv} \) is a set containing pairs of end-vertices of matching edges, and \( F^e \) is a set of faulty edges.  
Define \( F_i = (F^{mv} \cup F^e) \cap BP_n^i \) for \( i \in [n] \).  
Again, w.l.o.g., assume that \( |F_1| \geq |F_j| \) for all \( j \in [n] \setminus \{1\} \).  
We begin by considering the base condition.

\begin{lemma}\label{L13}
The $BP_{3}$ is $1$-hybrid fault Hamiltonian.
\end{lemma}

\begin{figure}[h]
  \centering
  \includegraphics[width=1\textwidth]{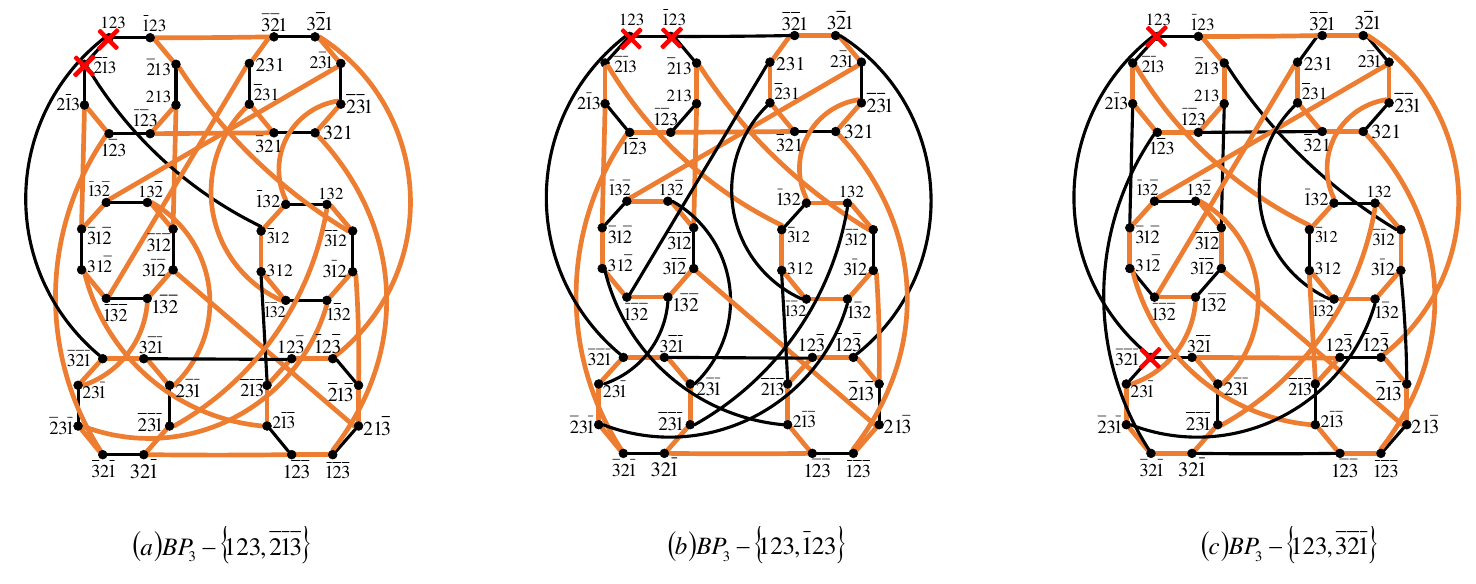}\\
  \caption{$BP_{3}-F^{mv}$ has a Hamiltonian cycle.}
  \label{Fig9}
\end{figure}

\noindent{\bf Proof.} In this case, we have \( |F| = |F^e| + |F^{mv}| \leq 1 \).  
If \( |F| = 0 \), then by Lemma \(\ref{L16}\), \( BP_3 \) contains a Hamiltonian cycle.     
Next, we consider the case where \( |F| = 1 \).  
If \( |F^e| = 1 \), then by Lemma \(\ref{L16}\), \( BP_3 \) still contains a Hamiltonian cycle.  
Otherwise, \( |F^{mv}| = 1 \).  

Since \( BP_3 \) is vertex-transitive, it suffices to analyze the case where the faulty vertex is \( 123 \). Based on its neighbors, we examine the following cases:  

\noindent{\bf Case 1.} The pair \( \{123, \bar{2}\bar{1}3\} \) is faulty.

The Hamiltonian cycle of $BP_{3}-\{123, \bar{2}\bar{1}3\}$ is $\langle \bar{1}23,\bar{3}\bar{2}1,231,\bar{1}\bar{3}\bar{2},31\bar{2},2\bar{1}\bar{3},\bar{2}\bar{1}\bar{3},12\bar{3},\bar{1}2\bar{3},3\bar{2}$\\
$1,2\bar{3}1,\bar{1}3\bar{2},\bar{3}1\bar{2},2\bar{1}3,1\bar{2}3,\bar{3}2\bar{1},\bar{2}3\bar{1},1\bar{3}2,3\bar{1}2,\bar{2}1\bar{3},21\bar{3},3\bar{1}\bar{2},1\bar{3}\bar{2},23\bar{1},\bar{3}\bar{2}\bar{1},3\bar{2}\bar{1},2\bar{3}\bar{1},13\bar{2},\bar{3}\bar{1}\bar{2},213,$\\
$\bar{1}\bar{2}3,\bar{3}21,\bar{2}31,\bar{1}\bar{3}2,312,\bar{3}12,\bar{1}32,\bar{2}\bar{3}1,321,\bar{1}\bar{2}\bar{3},1\bar{2}\bar{3},32\bar{1},\bar{2}\bar{3}\bar{1},132,\bar{3}\bar{1}2,\bar{2}13,\bar{1}23\rangle$(see Figure \ref{Fig9}($a$)).

\noindent{\bf Case 2.} The pair $\{123, \bar{1}23\}$ is faulty.

The Hamiltonian cycle of  $BP_{3}-\{123,\bar{1}23\}$ is $\langle
\bar{2}\bar{1}3,2\bar{1}3,\bar{3}1\bar{2},31\bar{2},\bar{1}\bar{3}\bar{2},1\bar{3}\bar{2},3\bar{1}\bar{2},21\bar{3},\bar{1}\bar{2}\bar{3},3$\\
$21,\bar{2}\bar{3}1,\bar{1}32,132,\bar{3}\bar{1}2,\bar{2}13,213,\bar{3}\bar{1}\bar{2},13\bar{2},\bar{1}3\bar{2},2\bar{3}1,3\bar{2}1,\bar{3}\bar{2}1,231,\bar{2}31,\bar{3}21,\bar{1}\bar{2}3,1\bar{2}3,\bar{3}2\bar{1},\bar{2}3\bar{1},23$\\
$\bar{1},\bar{3}\bar{2}\bar{1},3\bar{2}\bar{1},2\bar{3}\bar{1},\bar{2}\bar{3}\bar{1},32\bar{1},1\bar{2}\bar{3},2\bar{1}\bar{3},\bar{2}\bar{1}\bar{3},12\bar{3},\bar{1}2\bar{3},\bar{2}1\bar{3},3\bar{1}2,1\bar{3}2,\bar{1}\bar{3}2,312,\bar{3}12,\bar{2}\bar{1}3\rangle$ (see Figure \ref{Fig9}($b$)).

\noindent{\bf Case 3.} The pair $\{123, \bar{3}\bar{2}\bar{1}\}$ is faulty.

The Hamiltonian cycle of $BP_{3}-\{123, \bar{3}\bar{2}\bar{1}\}$ is $\langle \bar{2}\bar{1}3,2\bar{1}3,1\bar{2}3,\bar{1}\bar{2}3,213,\bar{2}13,\bar{1}23,\bar{3}\bar{2}1,3\bar{2}1,\bar{1}$\\
$2\bar{3},\bar{2}1\bar{3},21\bar{3},3\bar{1}\bar{2},\bar{3}\bar{1}\bar{2},13\bar{2},2\bar{3}\bar{1},3\bar{2}\bar{1},12\bar{3},\bar{2}\bar{1}\bar{3},312,\bar{1}\bar{3}2,1\bar{3}2,3\bar{1}2,\bar{3}\bar{1}2,132,\bar{2}\bar{3}\bar{1},32\bar{1},\bar{3}2\bar{1},\bar{2}3\bar{1},23$\\
$\bar{1},1\bar{3}\bar{2},\bar{1}\bar{3}\bar{2},231,\bar{2}31,\bar{3}21,321,\bar{1}\bar{2}\bar{3},1\bar{2}\bar{3},2\bar{1}\bar{3},31\bar{2},\bar{3}1\bar{2},\bar{1}3\bar{2},2\bar{3}1,\bar{2}\bar{3}1,\bar{1}32,\bar{3}12,\bar{2}\bar{1}3\rangle$ (see Figure \ref{Fig9}($c$)).
\hfill$\Box$\\

When constructing Hamiltonian cycles or paths in \( BP_n \), we can first focus on constructing Hamiltonian paths for the subgraph \( BP_n^I - F_I \), which is induced by \( \bigcup_{i \in I} (BP_n^i - F_i) \) for \( I \subseteq [n] \).  
Next, we consider the following lemmas.

\begin{lemma}\label{L17}
For \( n \geq 4 \), let \( I \subseteq [n] \) with \( |I| \geq 5 \) and 
\( j_1, j_2 \in I \), where \( j_1 \neq j_2 \).  
Let \( F^{mv} \) be a set containing pairs of end-vertices of matching edges in \( BP_n \), and let \( F^e \) be a faulty edge set in \( BP_n - F^{mv} \), where \( F = F^{mv} \cup F^e \) with \( |F| \leq n - 2 \).  
If \( u \in V(BP_n^{j_1}) \) and \( v \in V(BP_n^{j_2}) \), and \( BP_n^i - F_i \) is Hamiltonian connected for all \( i \in I \), then \( BP_n^I - F_I \) contains a Hamiltonian path \( P[u, v] \) between \( u \) and \( v \).
\end{lemma}

\noindent{\bf Proof.} Let \( I = \{k_1', k_2', \cdots, k_m'\} \) with \( m \geq 5 \).  
W.l.o.g., assume that \( j_1 = k_1' \) and \( j_2 = k_m' \).  

Applying Proposition \(\ref{P2}\), there exists a permutation \( k_1, k_2, \cdots, k_m \) of \( k_1', k_2', \cdots, k_m' \) such that \( k_1 = k_1' \), \( k_m = k_m' \), and \( k_j \neq \overline{k_{j+1}} \) for \( j \in \{1, 2, \cdots, m-1\} \).  

By Lemma \(\ref{L14}\), the size of the edge set \( E_{k_j, k_{j+1}}(BP_n) \) is \( |E_{k_j, k_{j+1}}(BP_n)| = (n-2)! \cdot 2^{n-2} \). Since  
\[
|E_{k_j, k_{j+1}}(BP_n)| - |F| \geq (n-2)! \cdot 2^{n-2} - (n-2) > 0 \quad \text{for } n \geq 4,
\]  
we can choose a fault-free edge \( (v_{k_j}, u_{k_{j+1}}) \) (\((v_{k_j}, u_{k_{j+1}}) \notin F^e\) and \(\{v_{k_j}, u_{k_{j+1}}\}\cap V(F^{mv})=\emptyset\)) between \( BP_n^{k_j} \) and \( BP_n^{k_{j+1}} \), where \( v_{k_j} \in V(BP_n^{k_j}) \) and \( u_{k_{j+1}} \in V(BP_n^{k_{j+1}}) \) for \( 1 \leq j \leq m-1 \).  

Since \( BP_n^i - F_i \) is Hamiltonian connected for \( i \in I \), there exists a Hamiltonian path \( P[u_{k_j}, v_{k_j}] \) between \( u_{k_j} \) and \( v_{k_j} \) in \( BP_n^{k_j} \) for \( 1 \leq j \leq m \).  

\begin{figure}[h]
  \centering
  \includegraphics[width=0.8\textwidth]{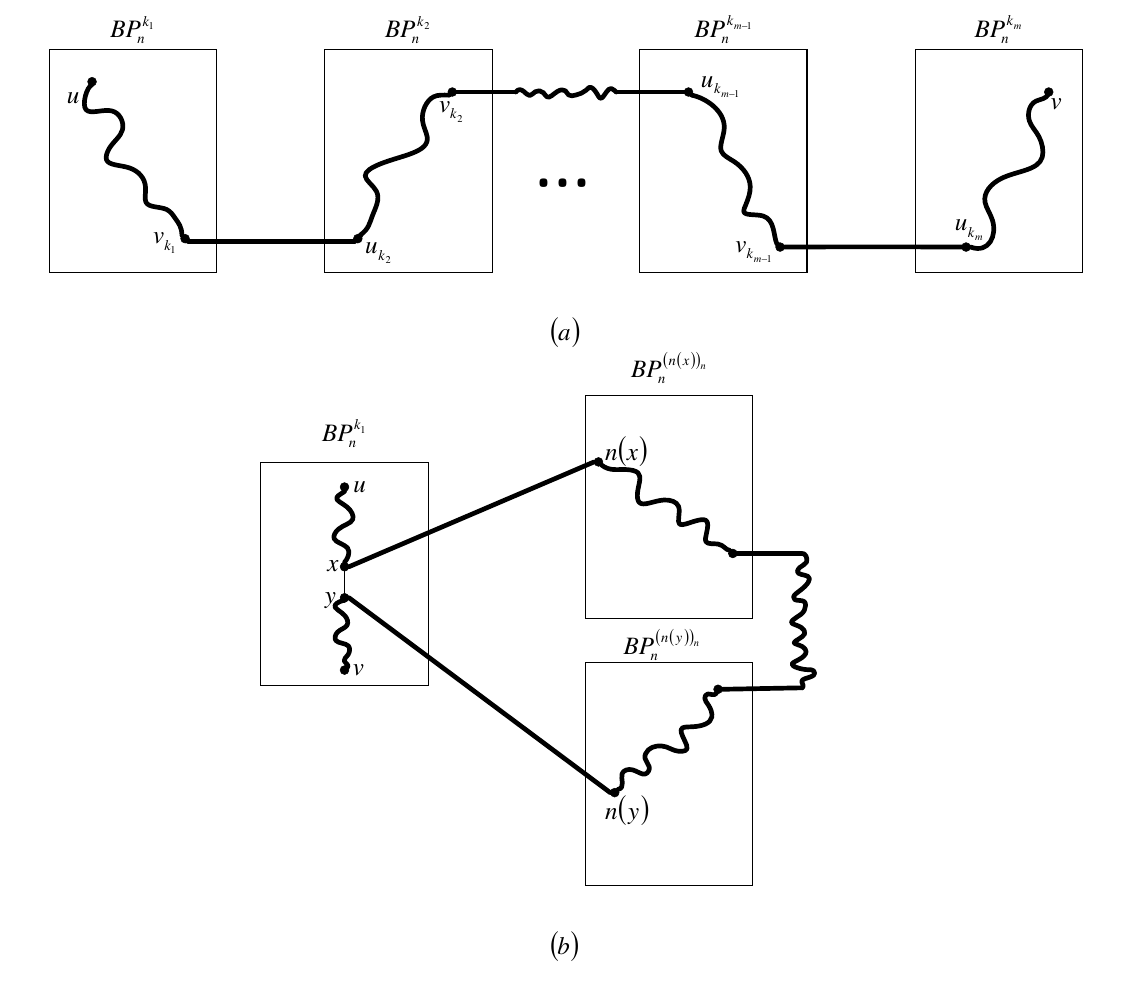}\\
  \caption{The illustration of Lemma \ref{L17} and Lemma \ref{L20}.}
  \label{Fig1}
\end{figure}

Thus, the Hamiltonian path of \( BP_n^I - F_I \) is given by  
\[
\langle u, P[u, v_{k_1}], v_{k_1}, u_{k_2}, P[u_{k_2}, v_{k_2}], \cdots, u_{k_m}, P[u_{k_m}, v], v \rangle,
\]  
as illustrated in Figure \ref{Fig1}$(a)$.

\hfill$\Box$\\

\begin{lemma}\label{L20}
For \( n \geq 4 \), let \( I \subseteq [n] \) with \( |I| \geq 6 \).  
Let \( F^{mv} \) be a set containing pairs of end-vertices of matching edges in \( BP_n \), and let \( F^e \) be a faulty edge set in \( BP_n - F^{mv} \). Let \( F = F^{mv} \cup F^e \) with \( |F| \leq n - 2 \).  
If \( u, v \in V(BP_n^{j_1}) \) with \( j_1 \in I \), and \( BP_n^i - F_i \) is Hamiltonian connected for \( i \in I \), then \( BP_n^I - F_I \) contains a Hamiltonian path \( P[u, v] \) between \( u \) and \( v \).  
\end{lemma}

\noindent{\bf Proof.}
Let \( I = \{k_1', k_2', \cdots, k_m'\} \) with \( m \geq 6 \).  
W.l.o.g., assume \( j_1 = k_1' \).  

Applying Proposition \(\ref{P2}\), there exists a permutation \( k_1, k_2, \cdots, k_m \) of \( k_1', k_2', \cdots, k_m' \) such that \( k_1 = k_1' \), \( k_m = k_m' \), and \( k_j \neq \overline{k_{j+1}} \) for all \( j \in \{1, 2, \cdots, m-1\} \).  

Within \( BP_n^{k_1} \), there are \( |V(BP_n^{k_1})| - |V(F^{mv}_{k_1})| \) fault-free vertices. Since \( BP_n^i - F_i \) is Hamiltonian connected for all \( i \in I \), there exists a Hamiltonian path \( P[u, v] \) in \( BP_n^{k_1} - F_{k_1} \), where the length of \( P[u, v] \) is \( |V(BP_n^{k_1})| - |V(F^{mv}_{k_1})| - 1 \).   

Outside \( BP_n^{k_1} \), there are \( |F| - |F_{k_1}| \) faulty elements. Since the out-neighbors of the two end-vertices associated with any element in \( F \setminus F_{k_1} \) belong to different two subgraphs, at most one out-neighbor of any element lies in \( BP_n^{k_1} \).  

For each element in \( F \setminus F_{k_1} \), each element is associated with at most one vertex of the Hamiltonian path in \( BP_n^{k_1} \) and with at most two edges of this path. Therefore, the number of edges in \( P[u, v] \) corresponding to \( V(F \setminus F_{k_1}) \) is at most \( 2(|F| - |F_{k_1}|) \). At the same time, the number of remaining edges in \( P[u, v] \) is at least \( |E(P[u, v])|-2(|F| - |F_{k_1}|) \). We claim that \( |E(P[u, v])|-2(|F| - |F_{k_1}|) >0\). This reasoning is as follows:

\text{For } $n \geq 4$,
\[
\begin{aligned}
|E(P[u, v])| - 2(|F| - |F_{k_1}|) 
&= (|V(BP_n^{k_1})| - |V(F^{mv}_{k_1})| - 1) - 2(|F| - |F_{k_1}|) \\
&= (|V(BP_n^{k_1})| - 1 - 2|F|) + (2|F_{k_1}|- |V(F^{mv}_{k_1})|) \\
&\geq |V(BP_n^{k_1})| - 1 - 2|F| \\
&\geq (n-1)! \cdot 2^{n-2} - 1 - 2(n-2) \\
&> 0 
\end{aligned}
\]
The first inequality follows from the definition of \( F \) and \( F^{mv} \), which implies \( |V(F^{mv}_{k_1})| \leq 2|F_{k_1}| \).

Thus, we can select a fault-free edge \( (x, y) \) in \( P[u, v] \) such that 
\( \{n(x), n(y)\} \cap V(F) = \emptyset \).  

Since \( (x)_n = (y)_n \) and \( d_{BP_n}(x, y) = 1 \), by Lemma \(\ref{L15}\), we have \( (n(x))_n \neq (n(y))_n \). Using Lemma \(\ref{L17}\), there exists a Hamiltonian path \( P[n(x), n(y)] \) in \( BP_n^{I'} - F_{I'} \), where \( I' = I \setminus \{k_1\} \).  

Thus, the Hamiltonian path of \( BP_n^I - F_I \) is  
\[
\langle u, P[u, x], x, n(x), P[n(x), n(y)], n(y), y, P[y, v], v \rangle,
\]  
as illustrated in Figure \ref{Fig1}$(b)$.
\hfill$\Box$\\

\begin{lemma}\label{L18}
For \( n \geq 4 \), if \( BP_{n-1} \) is \((n-4)\)-hybrid fault Hamiltonian connected and \((n-3)\)-hybrid fault Hamiltonian, then \( BP_n \) is \((n-2)\)-hybrid fault Hamiltonian.
\end{lemma}

\noindent{\bf Proof.} Let \( F^{mv} \) be a set containing pairs of end-vertices of matching edges in \( BP_n \), and let \( F^e \) be a set of faulty edges in \( BP_n - F^{mv} \), where \( F = F^{mv} \cup F^e \) and \( |F| \leq n-2 \).

If \( F^{mv} = \emptyset \), then by Lemma \(\ref{L16}\), the conclusion holds.

Next, we consider the case where \( F^{mv} \neq \emptyset \).  
Note that \( BP_n^i \) is isomorphic to \( BP_{n-1} \). By the symmetry of \( BP_n \), w.l.o.g., we assume that \( |F_1| \geq |F_j| \) for all \( j \in [n] \setminus \{1\} \).

Since \( |E_{i,j}(BP_n)| - |F| \geq (n-2)! \cdot 2^{n-2} - (n-2) > 0 \) for \( i, j \in [n] \) and \( i \notin \{j, \bar{j}\} \), there exists at least one fault-free edge between \( BP_n^i \) and \( BP_n^j \).  

\noindent{\bf Case 1.} \( |F_1| \leq n-4 \)

In this case, \( BP_n^i - F_i \) has a Hamiltonian path for \( i \in [n] \).  
We can choose a fault-free edge \( (u, n(u)) \) between \( BP_n^1 \) and \( BP_n^n \), where \( u \in V(BP_n^1) \) and \( n(u) \in V(BP_n^n) \).  
By Lemma \(\ref{L17}\), there exists a Hamiltonian path \( P[u, n(u)] \) in \( BP_n^I - F_I \) for \( I = [n] \).  
Thus, the Hamiltonian cycle of \( BP_n - F \) is  
\[
P[u, n(u)]\cup \{(n(u),u)\}.
\]

\noindent{\bf Case 2.} \( |F_1| = n-3 \)

Since \( |F_1| = n-3 \) and \( |F| \leq n-2 \), there exists at most one subgraph \( BP_n^i \) such that \( |F_i| \leq 1 \).  
By Lemma \(\ref{L14}\), there is no edge between \( BP_n^1 \) and \( BP_n^{\bar{1}} \).  
We now analyze the following subcases:

\begin{figure}[h]
  \centering
  \includegraphics[width=0.95\textwidth]{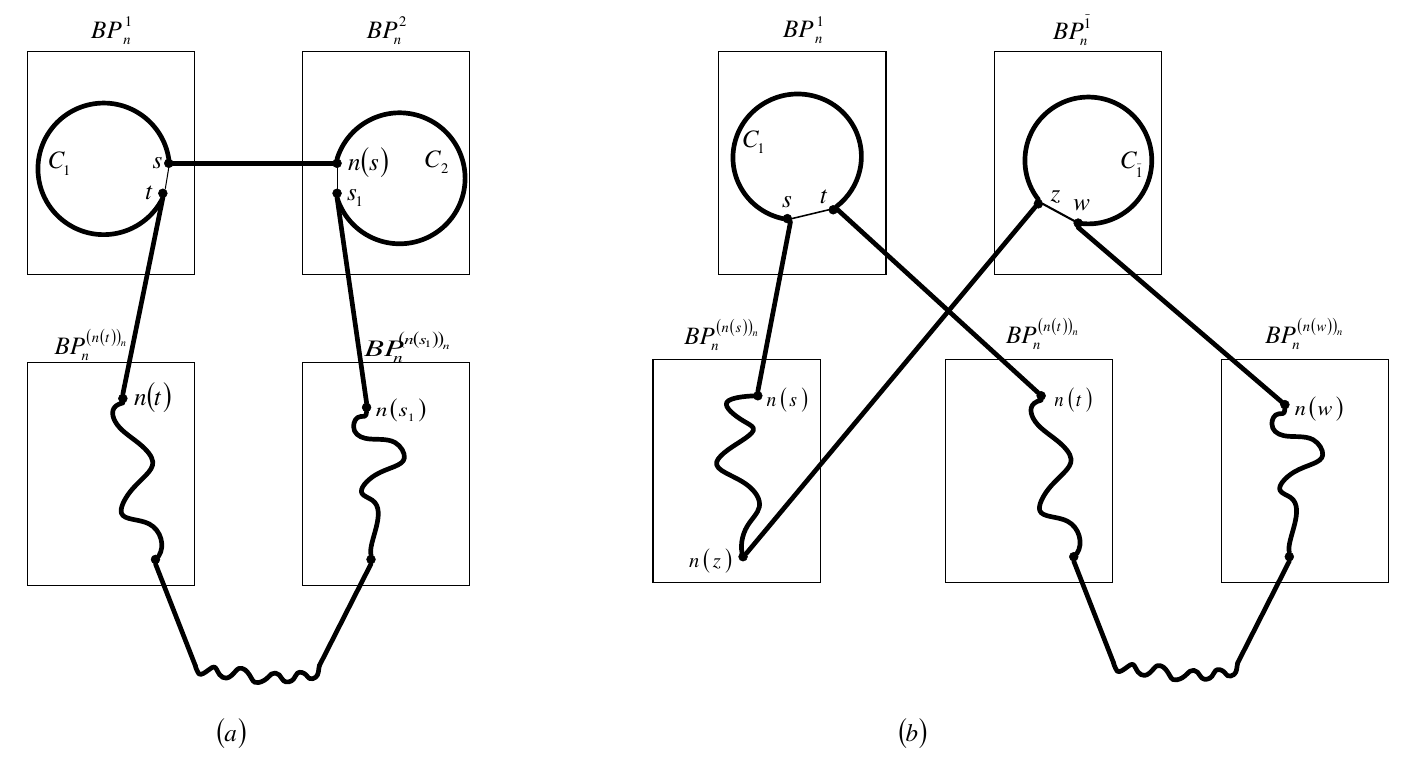}\\
  \caption{The illustration of \bf Case 2.}
  \label{Fig3}
\end{figure}

\noindent{\bf Subcase 2.1} \( i \neq \bar{1} \)

W.l.o.g., assume that \( i = 2 \), since \( |F_1| \geq |F_j| \) for \( j \in [n] \setminus \{1\} \).  
Since \( BP_{n-1} \) is \((n-3)\)-hybrid fault Hamiltonian, \( BP_n^1 - F_1 \) and \( BP_n^2 - F_2 \) have Hamiltonian cycles \( C_1 \) and \( C_2 \), respectively.  
Choose an edge \( (s, n(s)) \) between \( BP_n^1 \) and \( BP_n^2 \), where \( s \in V(C_1) \) and \( n(s) \in V(C_2) \).  

Additionally, let \( t \) be a neighbor of \( s \) on \( C_1 \), and \( s_1 \) be a neighbor of \( n(s) \) on \( C_2 \).  
Since \( 1 = (t)_n \neq (s_1)_n = 2 \) and \( d_{BP_n}(t, s_1) = 3 \), by Lemma \(\ref{L15}\), \( (n(s_1))_n \notin \{1, (n(t))_n\} \).  
Applying Lemma \(\ref{L17}\), there exists a Hamiltonian path \( P[n(t), n(s_1)] \) in \( BP_n^I - F_I \) for \( I = [n] \setminus \{1, 2\} \).  
Thus, the Hamiltonian cycle of \( BP_n - F \) is  
\[
\langle s, n(s), P[n(s), s_1], s_1, n(s_1), P[n(s_1), n(t)], n(t), t, P[t, s], s \rangle,
\]  
as illustrated in Figure \(\ref{Fig3}(a)\).

\noindent{\bf Subcase 2.2} \( i = \bar{1} \)

Since \( BP_{n-1} \) is \((n-3)\)-hybrid fault Hamiltonian, \( BP_n^1 - F_1 \) and \( BP_n^{\bar{1}} - F_{\bar{1}} \) have Hamiltonian cycles \( C_1 \) and \( C_{\bar{1}} \), respectively.  
Choose a vertex \( s \in V(C_1) \) and a vertex \( z \in V(C_{\bar{1}}) \) such that \( (n(z))_n = (n(s))_n \).  
Since every vertex in the cycle has degree two, choose a neighbor \( t \) of \( s \) on \( C_1 \) and a neighbor \( w \) of \( z \) on \( C_{\bar{1}} \), where \( (n(w))_n \neq (n(t))_n \).  
Since \( BP_{n-1} \) is \((n-4)\)-hybrid fault Hamiltonian connected, \( BP_n^{(n(s))_n} \) contains a Hamiltonian path \( P[n(s), n(z)] \).  
Applying Lemma \(\ref{L17}\), there exists a Hamiltonian path \( P[n(t), n(w)] \) in \( BP_n^I - F_I \) for \( I = [n] \setminus \{1, \bar{1}, (n(s))_n\} \). 
Thus, the Hamiltonian cycle of \( BP_n - F \) is  
\[
\langle s, n(s), P[n(s), n(z)], n(z), z, P[z, w], w, n(w), P[n(w), n(t)], n(t), t, P[t, s], s \rangle,
\]  
as illustrated in Figure \(\ref{Fig3}(b)\).

\noindent{\bf Case 3.} \( |F_1| = n-2 \)

Since \( F^{mv} \neq \emptyset \), we can select an element \( \{a_1, b_1\} \) from \( F^{mv}_1 \).  
As \( BP_{n-1} \) is \((n-3)\)-hybrid fault Hamiltonian, there exists a Hamiltonian cycle \( C_1 \) in \( BP_n^1 - F_1 \setminus \{\{a_1, b_1\}\} \).

\noindent{\bf Subcase 3.1.} $(a_{1},b_{1})\in E(C_{1})$.

In this case, \( C_1 \) can be expressed as  
\[
C_1 = \langle x_1, a_1, b_1, y_1, P[y_1, x_1], x_1 \rangle,
\]  
where \( x_1 \) is a neighbor of \( a_1 \) on \( C_1 \), and \( y_1 \) is a neighbor of \( b_1 \) on \( C_1 \).  

By Lemma \(\ref{L17}\), there exists a Hamiltonian path \( P[n(x_1), n(y_1)] \) in \( BP_n^I \) for \( I = [n] \setminus \{1\} \).  
Thus, the Hamiltonian cycle of \( BP_n - F^{mv} \) is  
\[
\langle x_1, n(x_1), P[n(x_1), n(y_1)], n(y_1), y_1, P[y_1, x_1], x_1 \rangle.
\]

\noindent{\bf Subcase 3.2.} $(a_{1},b_{1})\notin E(C_{1})$.

In this case, \( C_1 \) is given by  
\[
C_1 = \langle x_1, a_1, x_2, P[x_2, y_2], y_2, b_1, y_1, P[y_1, x_1], x_1 \rangle,
\]  
where \( x_1 \) and \( x_2 \) are neighbors of \( a_1 \) on \( C_1 \), \( y_1 \) and \( y_2 \) are neighbors of \( b_1 \) on \( C_1 \), and both \( x_1 \) and \( y_1 \) lie on the same path between \( a_1 \) and \( b_1 \) in \( C_1 \).  
Since \( (x_1)_n = (x_2)_n = 1 \) and \( d_{BP_n}(x_1, x_2) = 2 \), it follows from Lemma \(\ref{L15}\) that \( (n(x_1))_n \neq (n(x_2))_n \). Similarly, we have \( (n(y_1))_n \neq (n(y_2))_n \).  

Next, we analyze the following cases under the conditions \( (n(x_1))_n \neq (n(x_2))_n \) and \( (n(y_1))_n \neq (n(y_2))_n \).

\noindent{\bf Subcase 3.2.1.} $|\{(n(x_{1}))_{n},(n(x_{2}))_{n},(n(y_{1}))_{n},(n(y_{2}))_{n}\}|=4$.

\begin{figure}[h]
  \centering
  \includegraphics[width=0.9\textwidth]{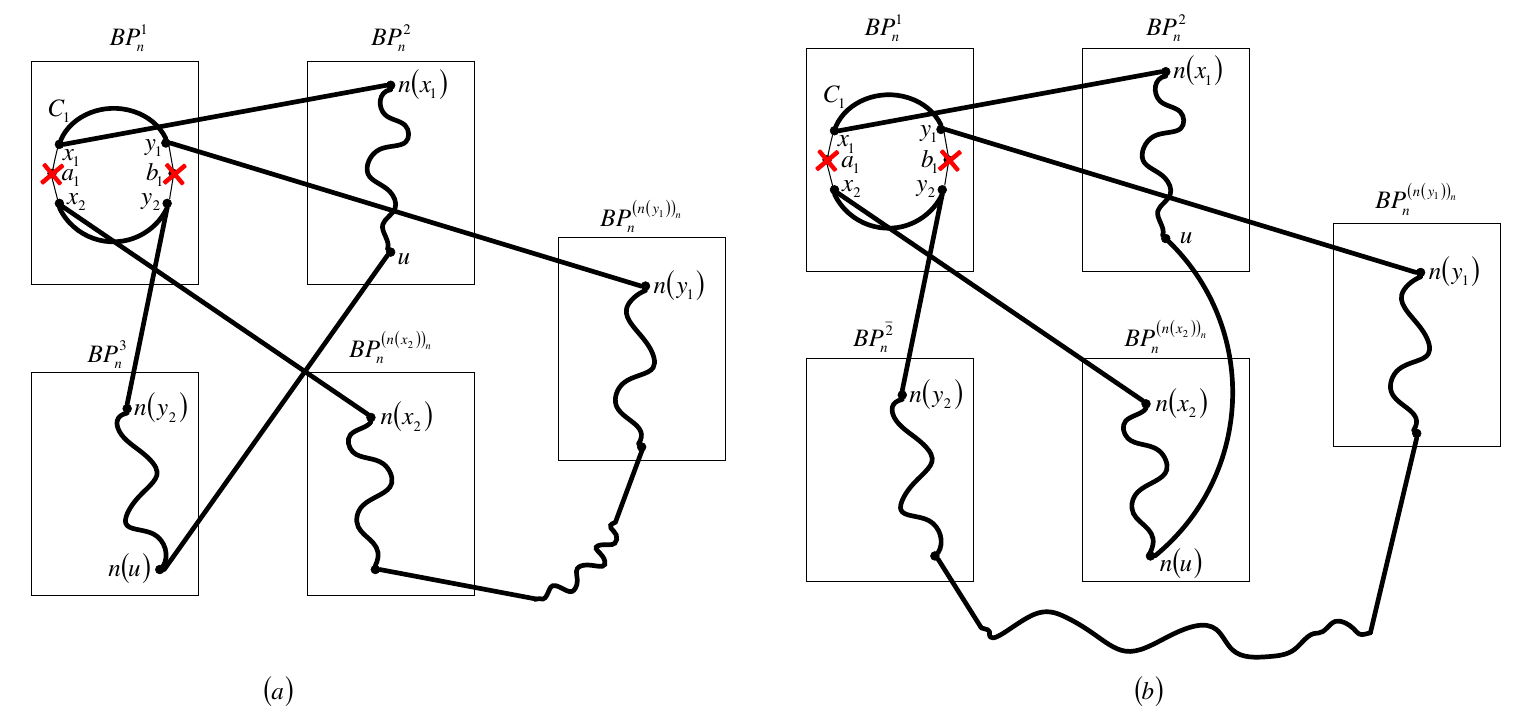}\\
  \caption{The illustration of \bf Subcase 3.2.1.}
  \label{Fig10}
\end{figure}

First, w.l.o.g, suppose that \( (n(x_1))_n = 2 \) and \( (n(y_2))_n = 3 \).  
We can select an edge \( (u, n(u)) \) between \( BP_n^2 \) and \( BP_n^3 \).  
Since \( BP_{n-1} \) is \((n-4)\)-hybrid fault Hamiltonian connected, there exists a Hamiltonian path \( P[n(x_1), u] \) in \( BP_n^2 \) and a Hamiltonian path \( P[n(u), n\)\\\((y_2)] \) in \( BP_n^3 \).  
Based on \( (n(x_2))_n \neq (n(y_1))_n \), by Lemma \(\ref{L17}\), we can find a Hamiltonian path \( P[n(x_2), n(y_1)] \) in \( BP_n^I \), where \( I = [n] \setminus \{1, 2, 3\} \).  
Thus, the Hamiltonian cycle of \( BP_n - F \) is  
\[
\begin{aligned}
\langle & x_1, n(x_1), P[n(x_1), u], u, n(u), P[n(u), n(y_2)], n(y_2), y_2, \\
& P[y_2, x_2], x_2, n(x_2), P[n(x_2), n(y_1)], n(y_1), y_1, P[y_1, x_1], x_1 \rangle.
\end{aligned}
\]
as illustrated in Figure \(\ref{Fig10}(a)\).

Now, consider the case where \( (n(x_1))_n = 2 \) and \( (n(y_2))_n = \bar{2} \).  
We can select an edge \( (u, n(u)) \) between \( BP_n^2 \) and \( BP_n^{(n(x_2))_n} \).  
Since \( BP_{n-1} \) is \((n-4)\)-hybrid fault Hamiltonian connected, there exists a Hamiltonian path \( P[n(x_1), u] \) in \( BP_n^2 \) and a Hamiltonian path \( P[n(u), n(x_2)] \) in \( BP_n^{(n(x_2))_n} \).  
By Lemma \(\ref{L17}\), there exists a Hamiltonian path \( P[n(y_2), n(y_1)] \) in \( BP_n^I \), where \( I = [n] \setminus \{1, 2, (n(x_2))_n\} \).  
Thus, the Hamiltonian cycle of \( BP_n - F \) is  
\[
\begin{aligned}
\langle & x_1, n(x_1), P[n(x_1), u], u, n(u), P[n(u), n(x_2)], n(x_2), x_2, \\
& P[x_2, y_2], y_2, n(y_2), P[n(y_2), n(y_1)], n(y_1), y_1, P[y_1, x_1], x_1 \rangle,
\end{aligned}
\]
as illustrated in Figure \(\ref{Fig10}(b)\).

\noindent{\bf Subcase 3.2.2.} $|\{(n(x_{1}))_{n},(n(x_{2}))_{n},(n(y_{1}))_{n},(n(y_{2}))_{n}\}|=3$.

\noindent{\bf Subcase 3.2.2.1.} $(n(x_{1}))_{n}=(n(y_{2}))_{n}$ or $(n(y_{1}))_{n}= (n(x_{2}))_{n}$. 

We provide the proof for the case where \( (n(x_1))_n = (n(y_2))_n \). The case \( (n(y_1))_n = (n(x_2))_n \) can be proven similarly.

W.l.o.g, suppose that \( (n(x_1))_n = (n(y_2))_n = 2 \). Since \( BP_{n-1} \) is \((n-4)\)-hybrid fault Hamiltonian connected, there exists a Hamiltonian path \( P[n(x_1), n(y_2)] \) in \( BP_n^2 \).  
By Lemma \(\ref{L17}\) and the condition \( (n(x_2))_n \neq (n(y_1))_n \), there exists a Hamiltonian path \( P[n(x_2), n(y_1)] \) in \( BP_n^I \), where \( I = [n] \setminus \{1, 2\} \).  
Thus, the Hamiltonian cycle of \( BP_n - F \) is  
\[
\begin{aligned}
\langle & x_1, n(x_1), P[n(x_1), n(y_2)], n(y_2), y_2, P[y_2, x_2], x_2, \\
& n(x_2), P[n(x_2), n(y_1)], n(y_1), y_1, P[y_1, x_1], x_1 \rangle.
\end{aligned}
\]


\noindent{\bf Subcase 3.2.2.2.} $(n(x_{1}))_{n}=(n(y_{1}))_{n}$ or $(n(x_{2}))_{n}=(n(y_{2}))_{n}$.

\begin{figure}[h]
  \centering
  \includegraphics[width=0.85\textwidth]{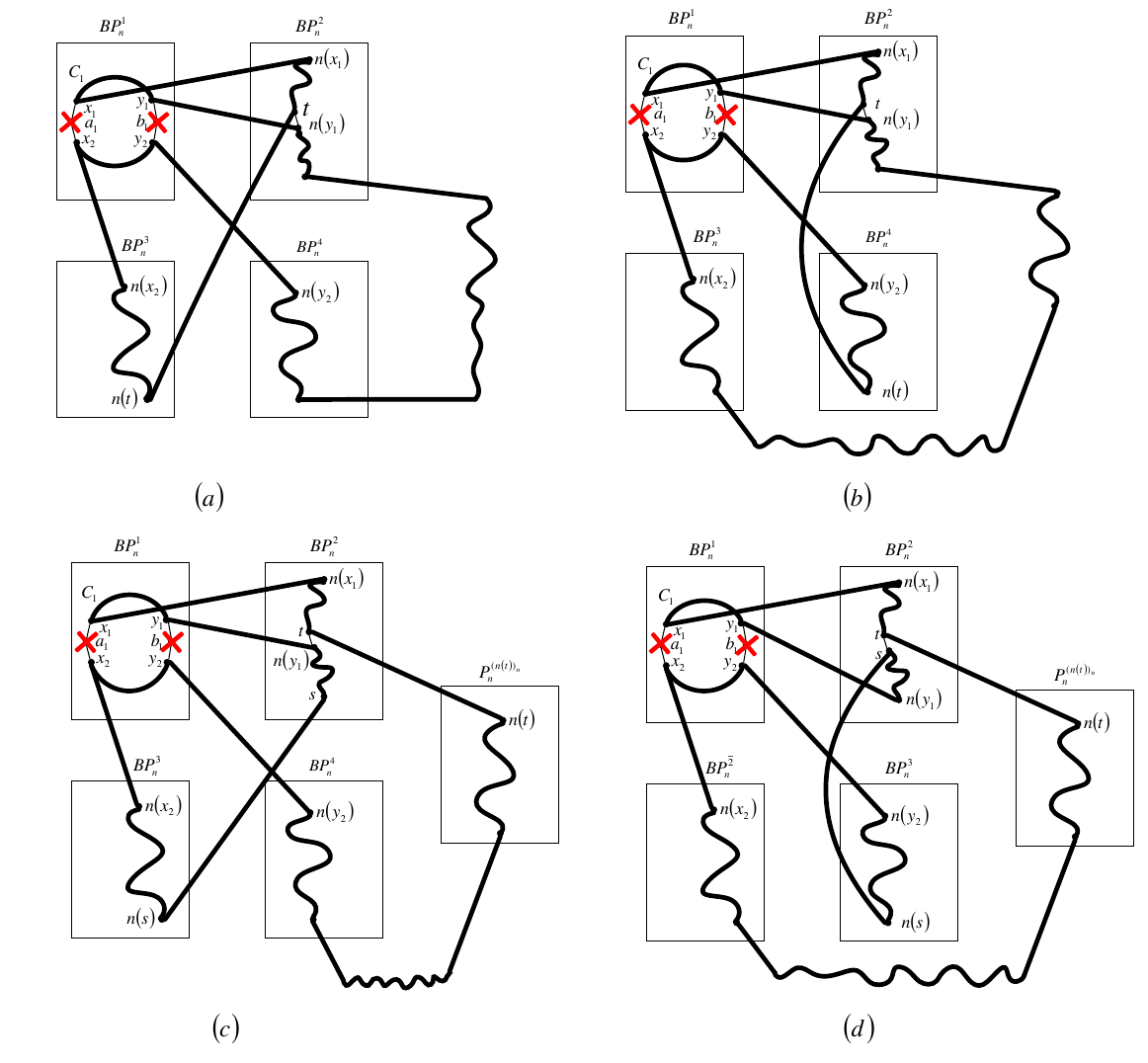}\\
  \caption{The illustration of \bf Subcase 3.2.2.2.}
  \label{Fig11}
\end{figure}
We provide the proof for the case where \( (n(x_1))_n = (n(y_1))_n\). The case \( (n(x_{2}))_{n}=(n(y_{2}))_{n} \) can be proven similarly.

First, w.l.o.g, suppose that \( (n(x_1))_n = (n(y_1))_n = 2 \), \( (n(x_2))_n = 3 \), and \( (n(y_2))_n = 4 \).  
Let \( t \) be a neighbor of \( n(y_1) \) on a Hamiltonian path in \( BP_n^2 \), such that \( t \) lies between \( n(x_1) \) and \( n(y_1) \) on the path, where the path starts with \( n(x_1) \) as its end-vertex.  
By Lemma \(\ref{L15}\), we have \( (n(t))_n \neq 1 \).  
Now, we consider the condition where \( (n(t))_n \in [n] \setminus \{1, \bar{2}\} \).

If \( (n(t))_n = 3 \), then by Lemma \(\ref{L17}\), there is a Hamiltonian path  
\[
P[n(x_1), n(y_2)] = \langle n(x_1), P[n(x_1), t], t, n(y_1), P[n(y_1), n(y_2)], n(y_2) \rangle
\]  
in \( BP_n^I \) for \( I = [n] \setminus \{1, 3\} \).  
Since \( BP_{n-1} \) is \((n-4)\)-hybrid fault Hamiltonian connected, there exists a Hamiltonian path \( P[n(x_2), n(t)] \) in \( BP_n^3 \).  
Thus, the Hamiltonian cycle of \( BP_n - F \) is  
\[
\begin{aligned}
 \langle &x_1, n(x_1), P[n(x_1), t], t, n(t), P[n(t), n(x_2)], n(x_2), x_2,P[x_2, y_2],\\ 
 & y_2, n(y_2), P[n(y_2), n(y_1)], n(y_1), y_1, P[y_1, x_1], x_1 \rangle   
\end{aligned}
\]  
as illustrated in Figure \(\ref{Fig11}(a)\).

If \( (n(t))_n = 4 \), then by Lemma \(\ref{L17}\), there is a Hamiltonian path  
\[
P[n(x_1), n(x_2)] = \langle n(x_1), P[n(x_1), t], t, n(y_1), P[n(y_1), n(x_2)], n(x_2) \rangle
\]  
in \( BP_n^I \) for \( I = [n] \setminus \{1, 4\} \).  
Since \( BP_{n-1} \) is \((n-4)\)-hybrid fault Hamiltonian connected, there exists a Hamiltonian path \( P[n(y_2), n(t)] \) in \( BP_n^4 \).  
Thus, the Hamiltonian cycle of \( BP_n - F \) is  
\[
\begin{aligned}
 \langle &x_1, n(x_1), P[n(x_1), t], t, n(t), P[n(t), n(y_2)], n(y_2), y_2, \\ &P[y_2, x_2], x_2, n(x_2), P[n(x_2), n(y_1)], n(y_1), y_1, P[y_1, x_1], x_1 \rangle   
\end{aligned}
\]  
as illustrated in Figure \(\ref{Fig11}(b)\).

If \( (n(t))_n \notin \{3, 4\} \), let \( s \) be a vertex in \( BP_n^2 \) such that \( (n(s))_n = 3 \).  
Since \( BP_{n-1} \) is \((n-4)\)-hybrid fault Hamiltonian connected, there exists a Hamiltonian path \( P[n(x_1), s] \) in \( BP_n^2 \) and a Hamiltonian path \( P[n(x_2), n(s)] \) in \( BP_n^3 \).  
By Lemma \(\ref{L17}\), there exists a Hamiltonian path \( P[n(y_2), n(t)] \) in \( BP_n^I \) for \( I = [n] \setminus \{1, 2, 3\} \).  
Thus, the Hamiltonian cycle of \( BP_n - F \) is  
\[
\begin{aligned}
\langle & x_1, n(x_1), P[n(x_1), t], t, n(t), P[n(t), n(y_2)], n(y_2), y_2, \\
& P[y_2, x_2], x_2, n(x_2), P[n(x_2), n(s)], n(s), s, P[s, n(y_1)], n(y_1), y_1, P[y_1, x_1], x_1 \rangle
\end{aligned}
\]  
as illustrated in Figure \(\ref{Fig11}(c)\).

Now, consider the case where \( (n(x_1))_n = (n(y_1))_n = 2 \) and \( \bar{2} \in \{(n(x_2))_n, (n(y_2))_n\} \).  
W.l.o.g., assume that \( (n(x_1))_n = (n(y_1))_n = 2 \), \( (n(x_2))_n = \bar{2} \), and \( (n(y_2))_n = 3 \).  
Since \( BP_{n-1} \) is \((n-4)\)-hybrid fault Hamiltonian connected, there exists a Hamiltonian path \( P[n(x_1), n(y_1)] \) in \( BP_n^2 \).  
As every vertex on the path has degree two except for the two end-vertices, we can select an edge \( (s, t) \) on \( P[n(x_1), n(y_1)] \) such that \( (n(s))_n = 3 \) and \( (n(t))_n \neq 1 \).  
Since \( BP_{n-1} \) is \((n-4)\)-hybrid fault Hamiltonian connected, there exists a Hamiltonian path \( P[n(y_2), n(s)] \) in \( BP_n^3 \).  
Given that \( (n(x_2))_n \neq (n(t))_n \) and applying Lemma \(\ref{L17}\), there exists a Hamiltonian path \( P[n(x_2), n(t)] \) in \( BP_n^I \), where \( I = [n] \setminus \{1, 2, \bar{2}\} \).  
Thus, the Hamiltonian cycle of \( BP_n - F \) is either  
\[
\begin{aligned}
\langle & x_1, n(x_1), P[n(x_1), t], t, n(t), P[n(t), n(x_2)], n(x_2), x_2, \\
& P[x_2, y_2], y_2, n(y_2), P[n(y_2), n(s)], n(s), s, P[s, n(y_1)], n(y_1), y_1, P[y_1, x_1], x_1 \rangle
\end{aligned}
\]  
or  
\[
\begin{aligned}
\langle & x_1, n(x_1), P[n(x_1), s], s, n(s), P[n(s), n(y_2)], n(y_2), y_2, \\
& P[y_2, x_2], x_2, n(x_2), P[n(x_2), n(t)], n(t), t, P[t, n(y_1)], n(y_1), y_1, P[y_1, x_1], x_1 \rangle.
\end{aligned}
\]

Figure \(\ref{Fig11}(d)\) illustrates the case where \( s \) is not on \( P[n(x_1), t] \).

\noindent{\bf Subcase 3.2.3.} $|\{(n(x_{1}))_{n},(n(x_{2}))_{n},(n(y_{1}))_{n},(n(y_{2}))_{n}\}|=2$.

\noindent{\bf Subcase 3.2.3.1.} $(n(x_{1}))_{n}=(n(y_{2}))_{n}$ and $(n(y_{1}))_{n}=(n(x_{2}))_{n}$.

W.l.o.g., assume that $(n(x_{1}))_{n}=(n(y_{2}))_{n}=2$.
Since $BP_{n-1}$ is $(n-4)$-hybrid fault Hamiltonian connected, there exists a Hamiltonian path $P[n(x_{1}),n(y_{2})]$ in $BP_{n}^{2}$.
Applying Lemma \ref{L20}, there exists a Hamiltonian path $P[n(x_{2}),n(y_{1})]$ in $P_{n}^{I}$, where $I=\langle n\rangle\setminus\{1,2\}$.
Thus, the Hamiltonian cycle of $BP_{n}-F$ is 
\[
\begin{aligned}
\langle &x_{1},n(x_{1}),P[n(x_{1}),n(y_{2})],n(y_{2}),y_{2},P[y_{2},x_{2}],x_{2},n(x_{2}),\\
&P[n(x_{2}),n(y_{1})],n(y_{1}),y_{1},P[y_{1},x_{1}],x_{1}\rangle.
\end{aligned}
\]
\noindent{\bf Subcase 3.2.3.2.}$(n(x_{1}))_{n}=(n(y_{1}))_{n}$ and $(n(x_{2}))_{n}=(n(y_{2}))_{n}$.

\begin{figure}[h]
  \centering
  \includegraphics[width=0.9\textwidth]{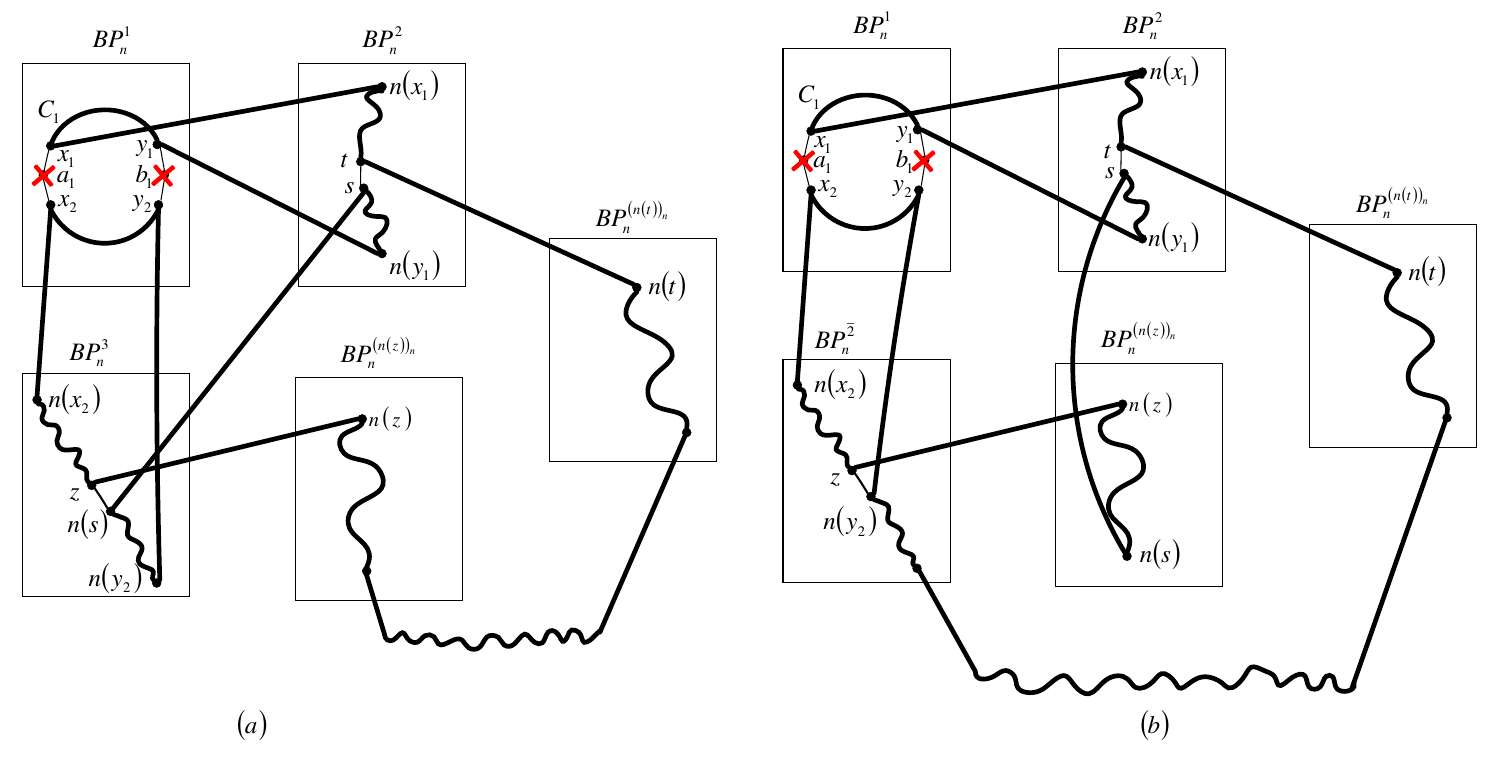}\\
  \caption{The illustration of case 3.2.3.2.}
  \label{F12}
\end{figure}

First, w.l.o.g., assume that $(n(x_{1}))_{n}=(n(y_{1}))_{n}=2$ and $(n(x_{2}))_{n}=(n(y_{2}))_{n}=3$.
Since $BP_{n-1}$ is $(n-4)$-hybrid fault Hamiltonian connected, there exists a Hamiltonian path $P[n(x_{1}),n(y_{1})]$ in $BP_{n}^{2}$ and $P[n(x_{2}),n(y_{2})]$ in $BP_{n}^{3}$, respectively.
We can select an edge \( (s, t) \) on \( P[n(x_1), n(y_1)] \) such that \( (n(s))_n = 3 \) and \( (n(t))_n \neq 1 \), as well as an edge \( (n(s), z) \) on \( P[n(x_2), n(y_2)] \) such that \( (n(z))_n \neq 1 \).  
Since \( d_{BP_n}(z, t) = 3 \) and \( 2 = (t)_n \neq (z)_n = 3 \), by Lemma \(\ref{L15}\), we conclude that \( (n(t))_n \neq (n(z))_n \).  
Applying Lemma \(\ref{L17}\), there exists a Hamiltonian path \( P[n(z), n(t)] \) in \( BP_n^I \), where \( I = [n] \setminus \{1, 2, 3\} \).  
Thus, the Hamiltonian cycle of \( BP_n - F \) is constructed as required (see Table \(\ref{B2}\)).  
Figure \(\ref{F12}(a)\) illustrates the case where \( s \) is not on \( P[n(x_1), t] \) and \( n(s) \) is not on \( P[n(x_2), z] \).

\begin{table}[h]
  \setlength{\abovecaptionskip}{0cm}
  \setlength{\belowcaptionskip}{0.2cm}
  \caption{the Hamiltonian cycle supposition 1 of Subcase 3.2.3.2 }\label{B2}
  \centering
  \begin{tabular}{>{\centering\arraybackslash}m{5cm} >{\raggedright\arraybackslash}m{11.8cm}}
    \hline
    the location of $s$ and $n(s)$ & the Hamiltonian cycle \\
    \hline
    $s$ is on $P[n(x_{1}),t]$ and $n(s)$ is on $P[n(x_{2}),z]$ &
    $\langle x_{1},n(x_{1}),P[n(x_{1}),s],s,n(s),P[n(s),n(x_{2})],n(x_{2}),x_{2},P[x_{2},y_{2}],y_{2},$\\
    &$n(y_{2}),P[n(y_{2}),z],z,n(z),P[n(z),n(t)],n(t),t,P[t,n(y_{1})],n(y_{1}),y_{1},$\\
    &$P[y_{1},x_{1}],x_{1}\rangle$ \\
    \hline
    $s$ isn't on $P[n(x_{1}),t]$ and $n(s)$ is on $P[n(x_{2}),z]$ &
    $\langle x_{1},n(x_{1}),P[n(x_{1}),t],t,n(t),P[n(t),n(z)],n(z),z,P[z,n(y_{2})],n(y_{2}),$ \\
    &$y_{2},P[y_{2},x_{2}],x_{2},n(x_{2}),P[n(x_{2}),n(s)],n(s),s,P[s,n(y_{1})],n(y_{1}),y_{1},$ \\
    &$P[y_{1},x_{1}],x_{1}\rangle$ \\
    \hline
    $s$ is on $P[n(x_{1}),t]$ and $n(s)$ isn't on $P[n(x_{2}),z]$ &
    $\langle x_{1},n(x_{1}),P[n(x_{1}),s],s,n(s),P[n(s),n(y_{2})],n(y_{2}),y_{2},P[y_{2},x_{2}],x_{2},$ \\
    &$n(x_{2}),P[n(x_{2}),z],z,n(z),P[n(z),n(t)],n(t),t,P[t,n(y_{1})],n(y_{1}),y_{1},$ \\
    &$P[y_{1},x_{1}],x_{1}\rangle$ \\
    \hline
    $s$ isn't on $P[n(x_{1}),t]$ and $n(s)$ isn't on $P[n(x_{2}),z]$ &
    $\langle x_{1},n(x_{1}),P[n(x_{1}),t],t,n(t),P[n(t),n(z)],n(z),z,P[z,n(x_{2})],n(x_{2}),$ \\
    &$x_{2},P[x_{2},y_{2}],y_{2},n(y_{2}),P[n(y_{2}),n(s)],n(s),s,P[s,n(y_{1})],n(y_{1}),y_{1},$ \\
    &$P[y_{1},x_{1}],x_{1}\rangle$ \\
    \hline
  \end{tabular}
\end{table}

Second, w.l.o.g., assume that $(n(x_{1}))_{n}=(n(y_{1}))_{n}=2$ and $(n(x_{2}))_{n}=(n(y_{2}))_{n}=\bar{2}$.
Let \( z \) be a neighbor of \( n(y_1) \) on a Hamiltonian path in \( BP_n^{\bar{2}} \) such that \( z \) lies between \( n(x_2) \) and \( n(y_2) \), where the path starts with \( n(x_2) \) as its end-vertex.  
Since \( BP_{n-1} \) is \((n-4)\)-hybrid fault Hamiltonian connected, there exists a Hamiltonian path \( P[n(x_1), n(y_1)] \) in \( BP_n^2 \).  
As every vertex on this path has degree two except the two end-vertices, we can select an edge \( (s, t) \) on \( P[n(x_1), n(y_1)] \) such that \( (n(s))_n = (n(z))_n \) and \( (n(t))_n \neq 1 \).  
By Lemma \(\ref{L17}\), there exists a Hamiltonian path  
\[
\langle n(x_2), P[n(x_2), z], z, n(y_2), P[n(y_2), n(t)], n(t) \rangle
\]  
in \( BP_n^I \), where \( I = [n] \setminus \{1, 2, (n(z))_n\} \).  
Since \( BP_{n-1} \) is \((n-4)\)-hybrid fault Hamiltonian connected, there exists a Hamiltonian path \( P[n(z), n(s)] \) in \( BP_n^{(n(z))_n} \).  
Thus, the Hamiltonian cycle of \( BP_n - F \) is either:  
\[
\begin{aligned}
\langle & x_1, n(x_1), P[n(x_1), t], t, n(t), P[n(t), n(y_2)], n(y_2), y_2, \\
& P[y_2, x_2], x_2, n(x_2), P[n(x_2), z], n(z), P[n(z), n(s)], n(s), s, P[s, n(y_1)], n(y_1), y_1, P[y_1, x_1], x_1 \rangle
\end{aligned}
\]  
or:  
\[
\begin{aligned}
\langle & x_1, n(x_1), P[n(x_1), s], s, n(s), P[n(s), n(z)], n(z), z, P[z, n(x_2)], \\
& n(x_2), x_2, P[x_2, y_2], y_2, n(y_2), P[n(y_2), n(t)], n(t), t, P[t, n(y_1)], n(y_1), y_1, P[y_1, x_1], x_1 \rangle.
\end{aligned}
\]  
Figure \(\ref{F12}(b)\) illustrates the case where \( s \) is not on \( P[n(x_1), t] \).

\hfill$\Box$\\

\begin{lemma}\label{L19}
For $n\geq 4$, if $BP_{n-1}$ is $(n-4)$-hybrid fault Hamiltonian connected and $(n-3)$-hybrid fault Hamiltonian, then $BP_{n}$ is $(n-3)$-hybrid fault Hamiltonian connected.
\end{lemma}

\noindent{\bf Proof.}  Let \( F^{mv} \) be a set containing pairs of end-vertices of matching edges in \( BP_n \), and let \( F^e \) be a set of faulty edges in \( BP_n - F^{mv} \), where \( F = F^{mv} \cup F^e \) and \( |F| \leq n-2 \).

If \( F^{mv} = \emptyset \), then by Lemma \(\ref{L16}\), the conclusion holds.

Next, we consider the case where \( F^{mv} \neq \emptyset \).

Since \( |E_{i,j}(BP_n)| - |F| \geq (n-2)! \cdot 2^{n-2} - (n-2) > 0 \) for \( i, j \in [n] \) and \( i \notin \{j, \bar{j}\} \), there is at least one edge between \( BP_n^i \) and \( BP_n^j \).  
Recall that \( BP_n^i \) is isomorphic to \( BP_{n-1} \).  
By the symmetry of \( BP_n \), w.l.o.g., we can assume that \( |F_1| \geq |F_j| \) for \( j \in [n] \setminus \{1\} \).  
Let \( u \in V(BP_n^{j_1}) \), \( v \in V(BP_n^{j_2}) \), and \( \{u, v\} \cap (V(F_{j_1}) \cup V(F_{j_2})) = \emptyset \), where \( j_1, j_2 \in [n] \).  
We will construct a Hamiltonian path between \( u \) and \( v \) in \( BP_n - F \).  
To proceed, we need to consider the following cases.

\noindent{\bf Case 1.} $|F_{1}|\leq n-4$.

In this case, $BP_{n}^{i}-F_{i}$ has a Hamiltonian path for all $ i\in [n]$.
By Lemmas \ref{L17} or \ref{L20}, there exists a Hamiltonian path $P[u,v]$ in $BP_{n}^{I}-F_{I}$, where $I=[n]$.

\noindent{\bf Case 2.} $|F_{1}|= n-3$.

Since $BP_{n-1}$ is $(n-3)$-hybrid fault Hamiltonian, there exists a Hamiltonian cycle $C_{1}$ in $BP_{n}^{1}-F_{1}$.

\noindent{\bf Subcase 2.1.} $|\{1,j_{1},j_{2}\}|=3$.

In this case, $j_{1}\neq j_{2}$, $j_{1}\neq1$ and $j_{2}\neq1$.

\begin{figure}[h]
  \centering
  \includegraphics[width=1\textwidth]{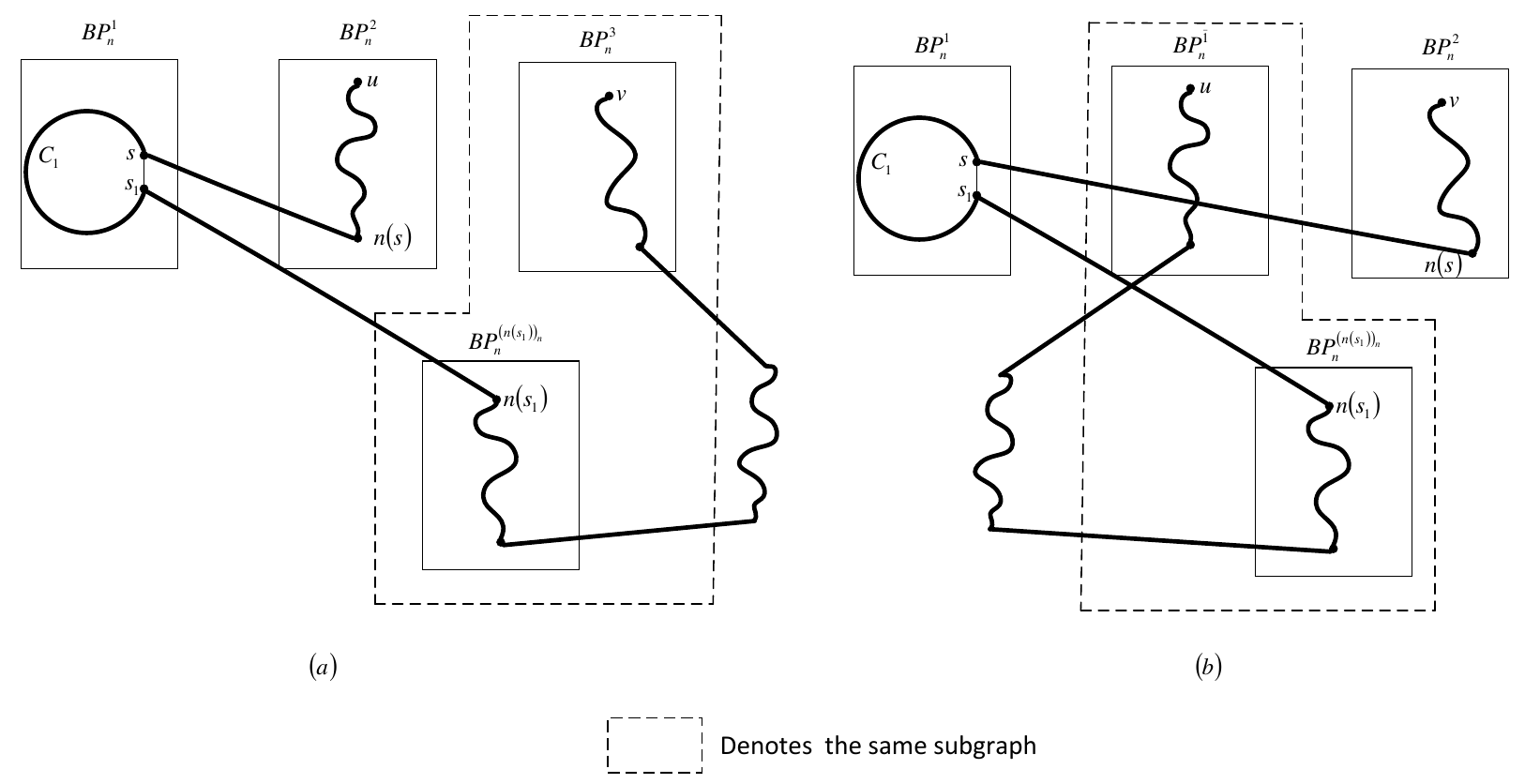}\\
  \caption{The illustration of case 2.1.}
  \label{Fig13}
\end{figure}

First, suppose that $\bar{1}\notin \{j_{1},j_{2}\}$.
W.l.o.g., assume that $j_{1}=2$ and $j_{2}=3$.
Choose a vertex \( s \) on \( C_1 \) such that \( (n(s))_n = 2 \), and let \( s_1 \) be a neighbor of \( s \) on \( C_1 \).  
By Lemma \(\ref{L15}\), \( (n(s_1))_n \neq (n(s))_n = 2 \).  
Since \( BP_{n-1} \) is \((n-4)\)-hybrid fault Hamiltonian connected, there exists a Hamiltonian path \( P[u, n(s)] \) in \( BP_n^2 \).  
By Lemmas \(\ref{L17}\) or \(\ref{L20}\), there exists a Hamiltonian path \( P[n(s_1), v] \) in \( BP_n^I \), where \( I = [n] \setminus \{1, 2\} \).  
Thus, the Hamiltonian path of \( BP_n - F \) is  
\[
\langle u, P[u, n(s)], n(s), s, P[s, s_1], s_1, n(s_1), P[n(s_1), v], v \rangle
\]  
as illustrated in Figure \(\ref{Fig13}(a)\).

Now, consider the case where \( \bar{1} \in \{j_1, j_2\} \).
W.l.o.g., assume that $j_{1}=\bar{1}$ and $j_{2}=2$.
Choose a vertex $s$ on $C_{1}$ such that $(n(s))_{n}=2$, and let $s_{1}$ be a neighbor of $s$ on $C_{1}$.
By Lemma \ref{L15}, $(n(s_{1}))_{n}\neq (n(s))_{n}=2$.
Since $BP_{n-1}$ is $(n-4)$-hybrid fault Hamiltonian connected, there exists a Hamiltonian path $P[v,n(s)]$ in $BP_{n}^{2}$.
By Lemmas \ref{L17} or \ref{L20}, there exists a Hamiltonian path $P[n(s_{1}),u]$ in $BP_{n}^{I}$, where $I=[n]\setminus \{1,2\}$.
Thus, the Hamiltonian path of $BP_{n}-F$ is 
\[\langle u,P[u,n(s_{1})],n(s_{1}),s_{1},P[s_{1},s],s,n(s),P[n(s),v],v\rangle\] 
as illustrated in Figure \(\ref{Fig13}(b)\).

\noindent{\bf Subcase 2.2. } $|\{1,j_{1},j_{2}\}|=2$.

\begin{figure}[h]
  \centering
  \includegraphics[width=1\textwidth]{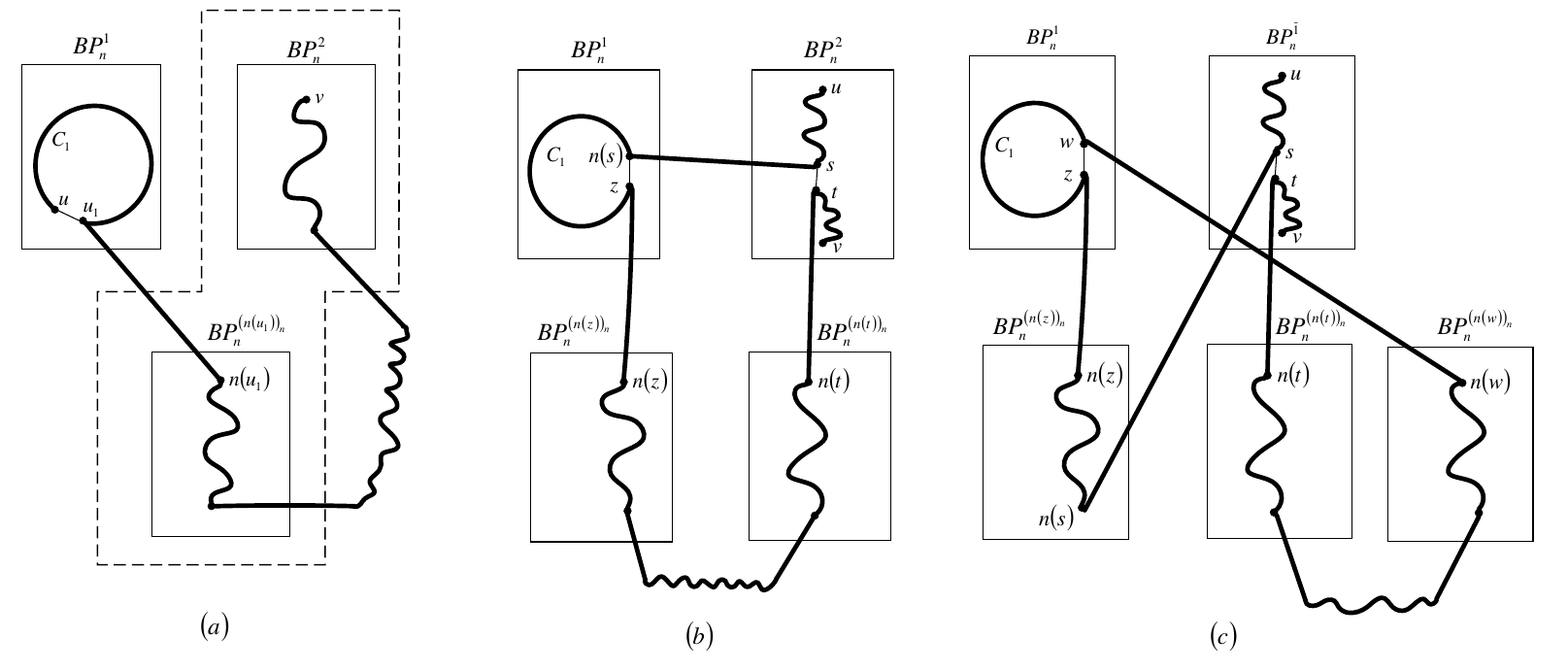}\\
  \caption{The illustration of case 2.2.}
  \label{F18}
\end{figure}

First, suppose that $1\in \{j_{1},j_{2}\}$.
W.l.o.g., assume that $j_{1}=1$ and $j_{2}=2$.
Choose a neighbor \( u_1 \) of \( u \) on \( C_1 \).  
By Lemmas \(\ref{L17}\) or \(\ref{L20}\), there exists a Hamiltonian path \( P[n(u_1), v] \) in \( BP_n^I \), where \( I = [n] \setminus \{1\} \).  
Thus, the Hamiltonian path of \( BP_n - F \) is  
\[
\langle u, P[u, u_1], u_1, n(u_1), P[n(u_1), v], v \rangle
\]  
as illustrated in Figure \(\ref{F18}(a)\).

Now, consider the case where $1\notin \{j_{1},j_{2}\}$.

If $j_{1}=j_{2}\in [n]\setminus\{\bar{1}\}$, w.l.o.g., assume that $j_{1}=j_{2}=2$.
Since $BP_{n-1}$ is $(n-4)$-hybrid fault Hamiltonian connected, there exists a Hamiltonian path $P[u,v]$ of $BP_{n}^{2}$.
We can select an edge \( (s, t) \) on \( P[u, v] \) such that \( (n(s))_n = 1 \).  
By Lemma \(\ref{L15}\), there is a neighbor \( z \) of \( n(s) \) on \( C_1 \) with \( (n(z))_n \neq 2 \).  
Since \( d_{BP_n}(z, t) = 3 \) and \( 1 = (t)_n \neq (z)_n = 2 \), by Lemma \(\ref{L15}\), we conclude that \( (n(t))_n \neq (n(z))_n \).  
Applying Lemma \(\ref{L17}\), there exists a Hamiltonian path \( P[n(z), n(t)] \) in \( BP_n^I \), where \( I = [n] \setminus \{1, 2\} \).  
Thus, the Hamiltonian path of \( BP_n - F \) is:  
\[
\langle u, P[u, s], s, n(s), P[n(s), z], z, n(z), P[n(z), n(t)], n(t), t, P[t, v], v \rangle
\]  
(see Figure \(\ref{F18}(b)\)).

If \( j_1 = j_2 = \bar{1} \), we need to construct a new Hamiltonian path for \( P[u, v] \).  
Since \( BP_{n-1} \) is \((n-4)\)-hybrid fault Hamiltonian connected, there exists a Hamiltonian path \( P[u, v] \) in \( BP_n^{\bar{1}} \).  
As every vertex on this path has degree two except the two end-vertices, we can select an edge \( (s, t) \) on \( P[u, v] \) and an edge \( (z, w) \) on \( C_1 \) such that \( (n(z))_n = (n(s))_n \) and \( (n(t))_n \neq (n(w))_n \). 
Since \( BP_{n-1} \) is \((n-4)\)-hybrid fault Hamiltonian connected, there exists a Hamiltonian path \( P[n(z), n(s)] \) in \( BP_n^{(n(s))_n} \).  
By Lemma \(\ref{L17}\), there exists a Hamiltonian path \( P[n(w), n(t)] \) in \( BP_n^I \), where \( I = [n] \setminus \{1, \bar{1}, (n(s))_n\} \).  
Thus, the Hamiltonian path of \( BP_n - F \) is either:  
\[
\langle u, P[u, s], s, n(s), P[n(s), n(z)], n(z), z, P[z, w], w, n(w), P[n(w), n(t)], n(t), t, P[t, v], v \rangle
\]  
or:  
\[
\langle u, P[u, t], t, n(t), P[n(t), n(w)], n(w), w, P[w, z], z, n(z), P[n(z), n(s)], n(s), s, P[s, v], v \rangle.
\]  
Figure \(\ref{F18}(c)\) illustrates the case where \( t \) is not on \( P[u, s] \).

\noindent{\bf Subcase 2.3.} $|\{1,j_{1},j_{2}\}|=1$.

In this case, $j_{1}=j_{2}=1$.

\noindent{\bf Subcase 2.3.1.} $(u,v)\in E(C_{1})$.

Choose an edge \( (s, t) \) on \( C_1 \) such that \( (s, t) \neq (u, v) \).  
By Lemma \(\ref{L17}\), there exists a Hamiltonian path \( P[n(s), n(t)] \) in \( BP_n^I \), where \( I = [n] \setminus \{1\} \).  
Thus, the Hamiltonian path of \( BP_n - F \) is:  
\[
\langle u, P[u, s], s, n(s), P[n(s), n(t)], n(t), t, P[t, v], v \rangle.
\]

\noindent{\bf Subcase 2.3.2.} $(u,v)\notin E(C_{1})$.

Choose a neighbor \( u_1 \) of \( u \) and a neighbor \( v_1 \) of \( v \) on \( C_1 \) such that \( u_1 \) and \( v_1 \) lie on different paths between \( u \) and \( v \) on \( C_1 \).  
By Lemmas \(\ref{L17}\) or \(\ref{L20}\), there exists a Hamiltonian path \( P[n(u_1), n(v_1)] \) in \( BP_n^I \), where \( I = [n] \setminus \{1\} \).  
Thus, the Hamiltonian path of \( BP_n - F \) is:  
\[
\langle u, P[u, v_1], v_1, n(v_1), P[n(v_1), n(u_1)], n(u_1), u_1, P[u_1, v], v \rangle.
\]
\hfill$\Box$\\

\begin{theorem}\label{T3}
For $n\geq 3$, $BP_{n}$ is $(n-2)$-hybrid fault Hamiltonian and $BP_{n}$ is $(n-3)$-hybrid fault Hamiltonian connected. Furthermore, all bounds are tight.
\end{theorem}

\noindent{\bf Proof.} We prove this theorem by induction on \( n \).  

{\bf Base case:} For \( n = 3 \), by Lemma \(\ref{L13}\), \( BP_3 \) is \( 1 \)-hybrid fault Hamiltonian.  
Additionally, by Lemma \(\ref{L16}\), \( BP_3 \) is Hamiltonian connected.  
Thus, the base case holds.  

{\bf Inductive step:} Assume the theorem holds for \( n-1 \).  
By Lemmas \(\ref{L18}\) and \(\ref{L19}\), we have proved the induction step.  

By the principle of mathematical induction, the theorem is true for \( n \geq 3 \).

Let \( F \) be the set of \( n-1 \) edges incident to a vertex of \( BP_n \). According to the \( n \)-regularity of \( BP_n \), removing these edges results in no Hamiltonian cycle in \( BP_n - F \). This demonstrates that the bound \( n-2 \) is tight.

Now, take \( F \) as the set of \( n-2 \) edges incident to a vertex \( u \) of \( BP_n \). Let the \( n-2 \) vertices (not $u$) corresponding to these edges be denoted as \( A \). Define \( \{x, y\} = N(u) \setminus A \), where \( x \) and \( y \) are the remaining neighbors of \( u \). In \( BP_n - F \), there is no Hamiltonian path between \( x \) and \( y \) since the unique path with $u$ as an inner vertex \( x \) and \( y \) as end-vertices is \( \langle x, u, y \rangle \). This demonstrates that the bound \( n-3 \) is tight.
\hfill$\Box$\\

\section{Conclusions}\label{section4}

In this paper, we have established the main result that \( BP_n \) is \((n-2)\)-hybrid fault Hamiltonian and \((n-3)\)-hybrid fault Hamiltonian connected for \( n \geq 3 \).  
In \cite{16}, Hung et al. demonstrated that the pancake graph \( P_n \) contains at least one Hamiltonian cycle (or path) under fault vertices and/or edges, respectively, where \( |F| \leq n-3 \) (or \( |F| \leq n-4 \)) and \( n \geq 4 \).  
In \cite{8}, Tsai et al. showed that pancake graph \( P_n \) contains at least one Hamiltonian cycle under the conditional fault edges, where \( |F| \leq 2n-7 \) and \( n \geq 4 \).  
However, to date, no work has been conducted on the Hamiltonian properties of burnt pancake graphs passing through prescribed edges under fault conditions.  
For future research, it would be intriguing to investigate the Hamiltonian properties of burnt pancake graphs that pass through prescribed edges while considering fault edges.

\section*{Acknowledgments}

The work was partially supported by
the Outstanding Youth Project of Hunan Provincial Department of Education(NO. 23B0117). This work is supported by NSFC (No. 12471326).

 \bibliographystyle{ieeetr}
\bibliography{main}

\end{document}